\documentclass[12pt]{article}
\usepackage[vmargin=1in,hmargin=1in]{geometry}

\usepackage{amssymb, cite}
\usepackage{graphicx}

\newcommand{\be}{\begin{equation}}
\newcommand{\ee}{\end{equation}}
\newcommand{\bqn}{\begin{eqnarray}}

\newcommand{\eqn}{\end{eqnarray}}

\newcommand{\bd}{\begin{description}}
\newcommand{\ed}{\end{description}}
\newtheorem{Theorem}{Theorem}[subsection]

\newtheorem{remark}[Theorem]{Remark}
\newtheorem{claim}[Theorem]{Claim}
\newtheorem{corollary}[Theorem]{Corollary}

\newtheorem{stat}{}[section]
\newtheorem{definition}[Theorem]{Definition}
\def\bs{\begin{stat}}
\def\es{\end{stat}}
\def\ben{\begin{enumerate}}
\def\een{\end{enumerate}}

\def\bp{\noindent{\bf Proof}  \ }
\newcommand{\ep}{\hfill $\square$}
\def\epcl{\hfill $\diamondsuit $}

\newcommand{\N}{{\mathbb N}}
\newcommand{\Z}{{\mathbb Z}}

\makeatletter
\@addtoreset{equation}{section}

\makeatother

\begin{document}

%\title{On k-Circular  Matroids of Graphs}

\begin{center}
{\large {\bf $K$-CIRCULAR MATROIDS OF GRAPHS}}
\\[4ex]
{\large {\bf Jos\'e F. De Jes\'us}}
\\[2ex]
{\bf University of Puerto Rico, San Juan, Puerto Rico, United States}
\\[4ex]
{\large {\bf Alexander Kelmans}}
\\[2ex]
{\bf University of Puerto Rico, San Juan, Puerto Rico, United States}
\\[0.5ex]
{\bf Rutgers University, New Brunswick, New Jersey, United States}

\end{center}

\date{}
\vskip 3ex

\begin{abstract}
 
 In 30's Hassler Whitney considered and completely solved the problem 
$(WP)$ of describing the classes of graphs $G$ having the same cycle matroid $M(G)$ \cite{W2,W3}. 
A natural analog $(WP)'$ of Whitney's problem $(WP)$ is  to describe the classes of graphs $G$ having the same  matroid $M'(G)$, where $M'(G)$ is a matroid (on the edge set of $G$) distinct from $M(G)$.
For example, the corresponding problem $(WP)'= (WP)_{\theta }$ for 
the so-called bicircular matroid $M_{\theta }(G)$ of graph $G$
was solved in \cite{CGW,Wag}.
  We define the so-called {\em $k$-circular matroid} $M_k(G)$ on the edge set of graph $G$ for any non-negative integer $k$ so that $M(G) = M_0(G)$ and $M_{\theta }(G) = M_1(G)$. It is natural to consider the corresponding analog $(WP)_k$ of Whitney's problem $(WP)$ not only for $k=0$ and $k=1$ but also for any integer $k \ge 2$. In this paper we give a characterization of the $k$-circular matroid $M_k(G)$ by describing the main constituents (circuits, bases, and cocircuits) in terms of graph $G$ and establish some important properties of the $k$-circular matroid. The results of this paper will be used in our further research on the problem $(WP)_k$. In our next paper we use these results to study a particular problem of $(WP)_k$ on graphs uniquely defined by their $k$-circular matroids.


\end{abstract}

\section{Introduction}

\indent

In 30's Hassler Whitney  developed a remarkable theory on the matroid isomorphism and the matroid duality of graphs \cite{W1,W2,W3,W4}.
He considered a graph $G$ and the so called cycle matroid $M(G)$ of $G$ 
(whose circuits  are the  edge subsets of the cycles in $G$)
and stated the following natural problems on pairs $\langle G, M(G)\rangle$:

\vskip 1ex

$(WP)$ describe the classes  of graphs having the same cycle matroid and, in particular, graphs that can be reconstructed
from cycle matroid (up to the names  of vertices) and

\vskip 1ex
$(WP^*)$
describe the pairs of graphs whose cycle matroids are dual, i.e. describe the class of graphs closed under their cycle matroids duality.

\vskip 1.5ex

Classical Whitney's graph matroid-isomorphism theorem and 
Whitney's planarity criterion provide the answers to the above questions \cite{W2,W3,W4} (see also \cite{Ox}).

Naturally, Whitney's problems and interesting results along this line prompted further questions and research on possible strengthenings as well as various extensions or analogs of some Whitney's results 
(see, for example, \cite{CGW,HalJ,HJK,Ksisd,K3skHng,Wag}).

\vskip 1.5ex

The goal of this research is twofold:
\\[1ex]
$(\Gamma 1)$
 to introduce and consider 
some new matroids related with a graph $G$ that are distinct from the cycle matroid $M(G)$ and to establish the structural properties of those matroids in terms of $G$ and
\\[1ex]
$(\Gamma 2)$ to extend some  of Whitney's results on problem 
$(WP)$
to this variety of graph matroids using  their structural properties established in $ (\Gamma 1)$

\vskip 1.5ex
This paper is the first one in the series of  our papers along the  $(\Gamma 1)$ - $(\Gamma 2)$ line on graphs and their matroids.

\vskip 1.5ex

Section  2 provides some basic notions, notation, and some known  facts on matroids and  graphs.

\vskip 1.5ex

In Section 3 some additional notions and some auxiliary and 
preliminary facts on graphs are presented that will be used later.
In particular, we will describe some useful properties of the 
function $\Delta (G) = |E(G)| - |V(G)|$ of graph $G$, establish
the splitter theorems for graphs with no leaves and no cycle components. We also define and characterize two special subgraphs of a graph that are called the {\em  core} and the {\em kernel}  of a graph.
These notions and results  will play an essential role in the study of the so-called $k$-circular matroids.

\vskip 1.5ex

In Section 4 we introduce the notion of the $k$-circular matroid $M_k(G)$ of a graph $G$, where $k$ is a non-negative integer, and establish some properties of this matroid.
In particular, we describe the main constituents of this matroid (bases, circuits, cocircuits, etc.) in graph terms.

\vskip 1.5ex
The results of this paper provide, in particular, a proper basis for our study of the problem 
$ (WP) _k$ on describing the classes of graphs with the same $k$-circular matroids.

\section{
Basic notions and facts on 
 matroids and graphs}

\label{basic} 

\subsection{On clutter and hereditary families}
\label{clutter}
\indent

Given a partial order set $(P, \preceq)$, 
a maximal element of $(P, \preceq)$ is called  {\em $\preceq $-maximal }  and 
a minimal element of $(P, \preceq)$ is called   {\em $\preceq $-minimal}.

The notions and facts described in this Section can be found in \cite{KmatR}. 
Given a finite set $E$, consider the poset ${\cal P} = (2^E, \subseteq)$  and 
${\cal X} \subseteq 2^E$. 

\vskip 1.5ex
Let 
\\
${\cal M}ax ({\cal X})$ denote the set of $\subseteq $-maximal elements of 
poset $({\cal X}, \subseteq )$ 
and

\vskip 0.7ex
\noindent
${\cal M}in ({\cal X})$ denote the set of $\subseteq $-minimal elements of 
poset $({\cal X}, \subseteq )$.

\vskip 1.5ex
Let
\\
${\cal M}ax ^{-1} ({\cal X}) =  \{X \subseteq E: X \subseteq Y~for~some ~Y \in 
~{\cal X}\}$
and
\vskip 0.7ex
\noindent
${\cal M}in ^{-1} ({\cal X}) =  \{X \subseteq E: X \supseteq Y~for~some ~Y \in 
~{\cal X}\}$.

\vskip 1.5ex
A family  ${\cal X}$ is called a {\em clutter} if 
$ {\cal M}ax ({\cal X}) = {\cal M}in ({\cal X})$, i.e. if $X \not \subseteq Y$ for every $X, Y  \in {\cal X}$ and $X \ne Y$.
Obviously, ${\cal M}ax ({\cal X})$ and ${\cal M}in ({\cal X})$ are clutters.

A family  ${\cal X}$ is called {\em hereditary } if 
$X \subseteq Y \in  {\cal X} \Rightarrow X \in  {\cal X}$.

A family  ${\cal X}$ is called {\em anti--hereditary } if 
$X \supseteq Y \in  {\cal X} \Rightarrow X \in  {\cal X}$.

\vskip 1.5ex
Obviously, the following is true.

\begin{claim}
 Let ${\cal X} \subseteq 2^E$. Then
\\[0.7ex]
$(c1)$ 
${\cal M}ax ^{-1} ({\cal X})$ is a hereditary family,
\\[0.7ex]
$(c2)$
${\cal M}in ^{-1} ({\cal X})$ is an anti-hereditary family,
\\[0.7ex]
$(c3)$
${\cal X}$ is a hereditary family if and only if 
${\cal M}in ^{-1}({\cal M}ax ({\cal X}) ) = {\cal X}$, and
\\[0.7ex]
$(c4)$
${\cal X}$ is an anti-hereditary family if and only if 
$ {\cal M}ax ({\cal M}in ^{-1} ({\cal X}) ) = {\cal X}$.
 
\end{claim}
 
Let ${\cal I} \subseteq 2^E$ and ${\cal I}$  a hereditary family and 
put $H = (E, {\cal I})$ and ${\cal I} = {\cal I} (H)$.
An element $X$ of $ {\cal I}$ is called an {\em independent set of $H$} and
${\cal I}$ is the {\em independence family of $H$}.

\vskip 1ex
Let ${\cal B} = {\cal B} (H) = {\cal M}ax ({\cal I})$.
An element $B$ of $ {\cal B}$ is called a {\em base of $H$} and
${\cal B} = {\cal B} (H)$ is the {\em  family of bases of $H$}.

\vskip 1ex
Let  ${\cal D}  = 2^E \setminus{\cal I}$.
Obviously,  $D \in {\cal D}$ if and only if $X$ is not independent set of $H$.
Therefore an element $D$ of ${\cal D}$ is called a  {\em dependent set  of $H$} and
${\cal D} = {\cal D} (H)$ is the {\em  dependence family  of $H$}.

Let ${\cal C}  = {\cal M}in ({\cal D)}$.
An element $C$ of $ {\cal C}$ is called a  {\em circuit  of $H$} and
${\cal C} = {\cal C} (H)$ is the {\em  family of circuits  of $H$}. A circuit $C$ of $H$ consisting of one element is called a {\em loop of $H$}.

\vskip 1ex
Given ${\cal B} = {\cal B} (H)$, let
 ${\cal B}^* = {\cal B}^*(H) = \{ E \setminus B: B \in {\cal B}\}$,
$ {\cal I}^* = {\cal M}ax ^{-1}({\cal B}^*)$,  
${\cal D}^*  = 2^E \setminus{\cal I}^*$, 
 ${\cal C}^*  = {\cal M}in ({\cal D}^*)$, and
 $H^* = (E, {\cal I}^*)$.
 We call $H = (E, {\cal I})$ and $H^* = (E, {\cal I}^*)$ {\em dual hereditary families on $E$} and put 
 ${\cal B} (H^*) =  {\cal B}^*(H)$,  ${\cal D} (H^*) =  {\cal D}^*(H)$, and
  ${\cal C} (H^*) =  {\cal C}^*(H)$. 

\vskip 1ex
 Accordingly, we call
 \\[0.4ex]
 an independent set of $H^*$ a {\em coindependet set of $H = (E, {\cal I})$},
\\[0.4ex]
a  dependent set $H^*$ a {\em codependent set  of $H = (E, {\cal I})$},
\\[0.4ex]
a base of $H^*$ a {\em cobase of $H = (E, {\cal I})$},
\\[0.4ex]
a circuit of  $H^*$ a {\em cocircuit of $H = (E, {\cal I})$}, and
\\[0.4ex]
a loop of $H^*$ a {\em coloop of $H = (E, {\cal I})$}.
 \\[1.4ex]
 \indent
 It is easy to see that the following is true.

\begin{claim}
\label{semi-matroid}

Let ${\cal I} \subseteq 2^E$, and ${\cal I}$ is a hereditary family.
Let  
$H = (E, {\cal I})$ and ${\cal I} = {\cal I} (H)$.
Then 
\\[0.7ex]
 $(c0)$ $  {\cal I} ^*(H)$  is a hereditary family,
\\[0.7ex]
 $(c1)$ $ {\cal B} (H)$ and ${\cal B}^*(H)$ are clutters,
 \\[0.7ex]
 $(c2)$  $ {\cal D} (H)$ and ${\cal D}^*(H)$  are  anti-hereditary families,
 \\[0.7ex]
  $(c3)$ $ {\cal C} (H)$ and  ${\cal C}^*(H)$ are clutters, and 
 \\[0.7ex] 
  $(c4)$  every family in  
  $\{ {\cal I}(H),~{\cal B}(H),~{\cal D} (H),~{\cal C} (H),~ 
   {\cal I}^*(H),~{\cal B}^*(H),~{\cal D}^* (H),~{\cal C}^* (H) \}$
   \\
  is uniquely  defined by any other family in the above list.
 
\end{claim}

 It is also easy to prove the following. 
 
\begin{claim}
\label{Base-Circuit}

 Let $H = (E, {\cal I} )$ be a non-empty hereditary family of subsets of $E$.
 Then 
\\[0.7ex]
 $(c1)$ $B \in  {\cal B} (H)~and~C \in  {\cal C} (H) \Rightarrow C \not \subseteq B$  %
 or, equivalentlly,
 \\[1ex]
 $(c2)$ $B \in  {\cal B} (H)~and~C^* \in  {\cal C} (H) \Rightarrow 
 B \cap  C^* \ne \emptyset $ and 
 \\
 $~~~~~~~B^* \in  {\cal B}^* (H)~and~C \in  {\cal C} (H) \Rightarrow  
  B^*\cap  C \ne \emptyset $.
 
\end{claim}

\subsection{On matroids}
\label{Matroids}

\indent

In this Section we describe some notions and necessary facts on matroids. Most of them can be found in \cite{Ox,Welsh}. Given a set $X \subseteq E$ and $e \in E$,  we usually write $X \cup e$ instead of $X \cup \{e\}$ and $X \setminus e$ instead of 
$X \setminus \{e\}$.

\vskip 1ex

A {\em matroid} is a pair $M = (E,{\cal I})$, where $E$ is a finite non-empty set and ${\mathcal I}\subseteq 2^{E}$ such that  
\\[0.7ex]
  $(A{\cal I}1)$  $  \emptyset \in {\cal I}$,  
 \\[0.7ex]
 $(A{\cal I}2)$  $X \subseteq Y \in {\cal I} \Rightarrow X \in {\cal I}$, and   
 \\[0.7ex]
$(A{\cal I}3)$  $X, Y \in {\cal I}~and~|X| < |Y| \Rightarrow \exists ~y \in Y \setminus X~ such ~that ~X \cup y \in {\cal I}$. 
  \\[0.7ex]
 A set $I \in {\cal I}$ is called an {\em  independent set of} $M$.
 
 \vskip 1ex
 Notice that   $(A{\cal I}1)$ and $(A{\cal I}2)$ imply that ${\cal I}$ is a hereditary family (see Section \ref{clutter}).
 Therefore $M = (E,{\cal I})$ is a matroid  if and only if 
 ${\cal I}$ is a non-empty hereditary family satisfying $(A{\cal I}3)$.
 Thus, all notions and claims about  hereditary families in Section \ref{clutter} are valid for the matroids.

\vskip 1ex
The next two claims provide criteria  for a family ${\cal F} \subseteq 2^E$  to be the family of bases and the family of circuits of a matroid, respectively.
\begin{claim}
\label{BasesOfM}
 
 Let ${\cal B} \subseteq 2^E$. Then 
 ${\cal B}$ is the family of bases of a matroid if and only if
 \\[0.7ex]
 $(A{\cal B}1)$ ${\cal B}$ is a non-empty clutter and
  \\[0.7ex]
$(A{\cal B}2)$ $X, Y \in {\cal B} ~and ~ X \ne Y \Rightarrow \forall ~
x \in X \setminus Y ~\exists ~ y \in Y \setminus X ~such~that~
 (X \setminus x) \cup y \in  {\cal B}$. 
 
\end{claim}  

From Claim \ref{BasesOfM} it follows that 
every two bases of matroid $M$  are of the same size.
Put $\rho (M)  =  |B|$ for $B \in {\cal B}(M)$ and 
$\rho ^*(M)  = |E| - \rho (M)$. 
We call $\rho (M)$  the {\em rank of matroid $M$} and 
$\rho ^*(M)$  the {\em corank of matroid $M$}.

\begin{figure}[h]
\begin{center}
\scalebox{0.3}[.3]{\includegraphics{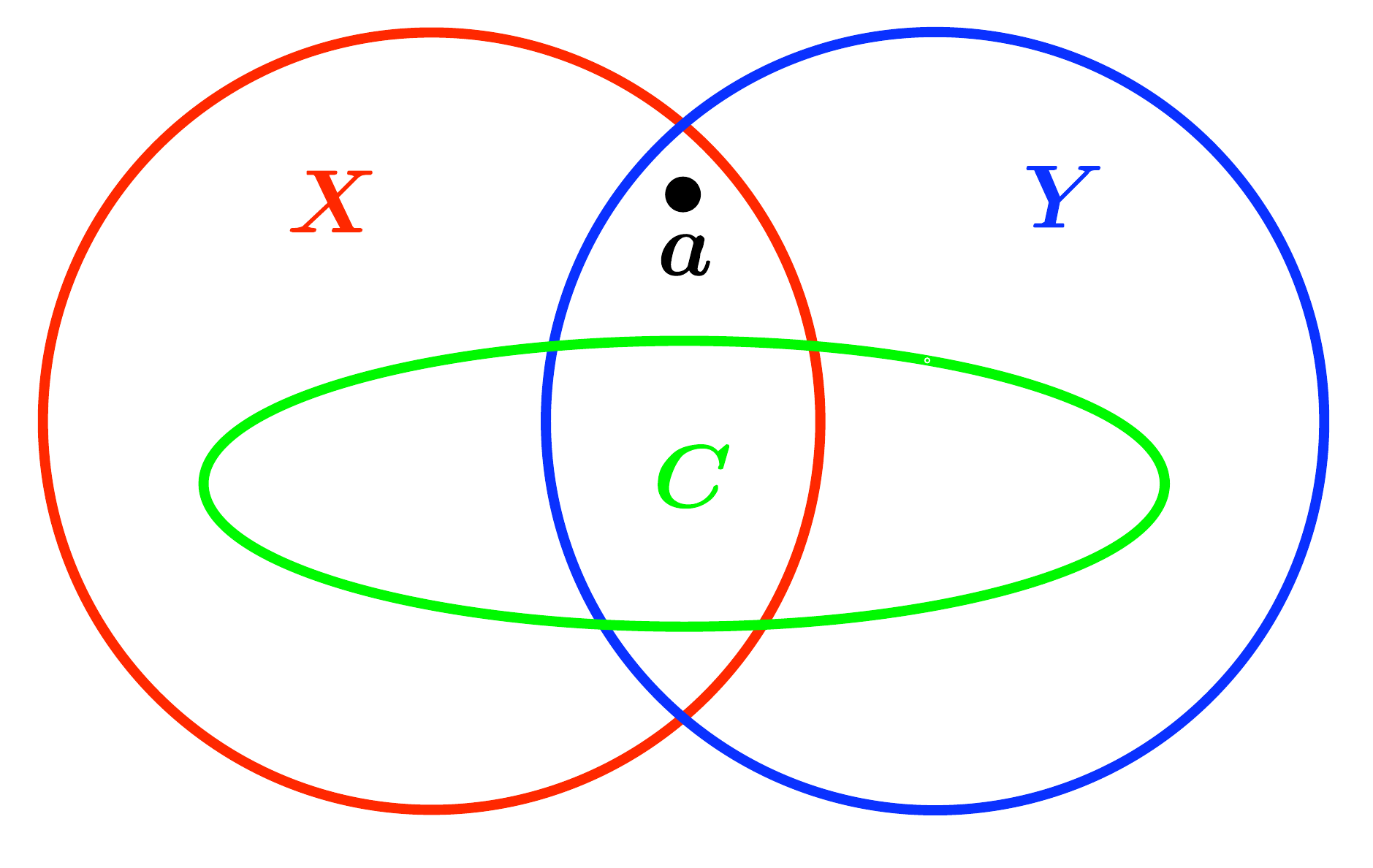}}
\end{center}
\caption
{Axiom $(A{\cal C}2)$ of a matroid}
\label{TsPresCirc}
\end{figure}

\begin{claim}

 Let ${\cal C} \subseteq 2^E$. Then 
 ${\cal C}$ is the family of circuits of a matroid if and only if
 \\[0.7ex]
 $(A{\cal C}1)$ ${\cal C}$ is a  clutter such that $\emptyset \not \in  {\cal C}$ and
  \\[0.7ex]
$(A{\cal C}2)$ $X, Y \in {\cal C} ~and ~ X \cap  Y \ne \emptyset \Rightarrow \forall ~
a \in X \cap Y ~\exists ~ C \in  {\cal C}~such~that~
 C \subseteq  (X \cup Y) \setminus a$ 
 \\
 $(${\em see Figure \ref{TsPresCirc}}$)$.

\end{claim}

 We will need the following simple  facts.

\begin{claim}
\label{MatroidColoops}  

 Let $M =  (E,{\cal I})$ be a matroid.
Then the following are equivalent: 
\\[1ex]
$(a1)$ $c$ is a coloop of $M$,
\\[1ex]
$(a2)$ $c \in B$ for every $B \in {\cal B}(M)$, and
\\[1ex]
$(a3)$ $c \not \in C$ for every $C \in {\cal C}(M)$.
 
\end{claim}

\begin{claim}
\label{FundCircuit}

  Let $M =  (E,{\cal I})$ be a matroid.
Then
\\[0.7ex]
$(c1)$ if $B  \in {\cal B}(M)$ and $e \in E \setminus B = B^*$, then
there exists a unique circuit $C  = C(e,B)$ of $M$ such that $e \in C \subseteq B \cup e$ $($or, equivalently, such that $C \cap B^* = \{e\}$$)$ and similarly,
 if $B  \in {\cal B}(M)$ and $e \in B$, then
there exists a unique cocircuit $C^*  = C^*(e,B)$ of $M$ such that $e \in C^* \subseteq B^* \cup e$ $($or, equivalently, such that $C^* \cap B = \{e\}$$)$,
$(c2)$
$u \in C(e,B) \setminus e \Leftrightarrow 
(B \setminus u) \cup e \in  {\cal B}(M)$ and similarly,
$u \in C^*(e,B) \setminus e \Leftrightarrow 
(B \setminus e) \cup u \in  {\cal B}(M)$,
\\[1ex]
$(c3)$ for every $C \in {\cal C} (M)$ $($$C^* \in {\cal C} (M)$$)$ there exists 
$B \in {\cal B}$ and $e \in E \setminus B = B^*$ such that $C  = C(e,B)$
$($respectively, 
$e \in  B$ such that $C^*  = C^*(e,B)$$)$.
 
\end{claim}

As in Section \ref{clutter}, given a matroid $M = (E, {\cal I})$, let
$M^* = (E, {\cal I}^*)$ be the pair with   ${\cal B}^*(M) = \{E \setminus B: B \in B(M)\}$.
It is easy to see that $M^*$  is a matroid and $M^*$ is the dual of $M$.

\vskip 1ex
Given  $B  \in {\cal B}(M)$ and $e \in E \setminus B = B^*$, we call 
$C(e,B)$ the {\em fundamental circuit of $B$  in $M$ rooted at $e$} or simply, the {\em $(B,e)$-circuit in $M$}.

\vskip 1ex
Similarly, if    $B^*  \in {\cal B}(M^*)$ and 
$e \in E \setminus B^* = B $, we call 
$C(e,B^*)$ the {\em fundamental circuit of $B^*$ in $M^*$ rooted at $e$} or simply, the {\em $(B^*,e)$-circuit in $M^*$}.

\vskip 1ex
We also call
$C(e,B^*)$ the {\em fundamental cocircuit of base $B$ in  $M$ rooted at $e$}, or simply, the {\em $(B,e)$-cocircuit in $M$}. 

\vskip 1ex

\begin{claim}
\label{EquivRel}

Let $M =  (E,{\cal I})$  be a matroid and $a, b \in E $, $a \ne b$.
 Then the following are equivalent:
 \\[1ex]
 $(a1)$  $a$ and $b$ belong to a circuit of $M$ and 
 \\[1ex]
 $(a2)$ $a$ and $b$ belong to a cocircuit of $M$.
 
\end{claim}  

Let $M =  (E,{\cal I})$  be a matroid, $L$ the set of loops, and  $L ^*$ the set of coloops of $M$. 
Given elements $a$ and  $b$  in 
$E \setminus (L \cup L ^*)$ we write $a \sim b$ if $a$ and $b$ belong to a common circuit of $M$ or, equivalently (by the above Claim),  if $a$ and $b$ belong to a common cocircuit of $M$.

\begin{claim}
\label{EquivRel}

$(E \setminus (L \cup L ^*), \sim )$ is an equivalence relation.
 
\end{claim}

\begin{claim}
\label{MatroidRestriction}  

 Let $M =  (E,{\cal I})$ be a matroid, $Z \subseteq E$, and  $M|_Z = (Z, \{X \subseteq Z: X \in {\cal I}\})$. 
Then $M|_Z $ is a matroid.
 
\end{claim}

For $X \subseteq E$   let 
 $M \setminus X = M|_{E\setminus X}$ and
$M / X =  (M^* \setminus X)^*$. 

\vskip 1.5ex 

A matroid $N$ is called a {\em component of matroid $M = (E, {\cal I})$}  if
$N = M|_Z$ for some equivalence class $Z$ of the equivalence relation 
$(E  \setminus (L \cup L ^*), \sim )$,
of matroid $M$, and so $N$ has at least two elements. Let $Cmp (M)$ denote the set of components of $M$.

\vskip 1.5ex
 A matroid $M = (E, {\cal I})$ is called {\em connected} 
 if $|E| \ge 2$ and $a \sim b$ for every $a, b \in E$, i.e. matroid $M = (E, {\cal I})$ is connected if and only if $M$ has exactly one components and $M$ has no loops and no coloops.

\begin{claim}
\label{Mcomponents}

Matroid $N$ is a component of matroid $M$
if and only if $N$ is a component of 
$M^*$, and so 
$Cmp (M) = Cmp (M^*)$.

\end{claim}

\vskip 1ex

We call a circuit $C$ of a connected matroid $M$ a {\em non-separating circuit } of $M$ if $M /C$ is a connected matroid.
Similarly, we call a cocircuit $C^*$  of a connected matroid $M$ a {\em non-separating cocircuit } of $M$
 if $M \setminus C^*$ is a connected matroid or, equivalently, if $C^*$ is a non-separating circuit of $M^*$.

\vskip 1.5ex

We call a matroid $M = (E, {\cal I})$ {\em non-trivial} if $E$ is neither a base of $M$ and nor a base of $M^*$ and {\em trivial}, otherwise. Obviously, $M$ is non-trivial if and only if ${\cal C}(M) \ne \emptyset$ and ${\cal C}^*(M) \ne \emptyset$. 

\vskip 1.5ex

Let $M = (E, {\cal I})$ and $M' = (E', {\cal I}')$ be matroids.
An {\em isomorphism from $M$ to $M'$ } is a bijection $\varepsilon  $ from $E$ to $E'$ such that $A \in {\cal I} \Leftrightarrow \varepsilon [A] \in {\cal I}$, where
$ \varepsilon [A] = \{ \varepsilon (e): e \in A\}$.
{\em Matroids $M$ and $M'$ are isomorphic} if there exists an isomorphism from $M$ to $M'$ (or, equivalently, an  isomorphism from $M'$ to $M$).

\subsection{On graphs}
\label{OnGraphs}

\indent

In this Section we describe some notions and necessary facts on graphs.
Most of them can be found in  \cite{BM,D,Wst}.

 \vskip 1.5ex
 
A {\em graph} $G$ is a triple $(V, E,\phi)$ such that $V$ and $E$ are disjoint finite sets, 
$V\cap E = \emptyset $, $V \ne \emptyset $, and $\phi: E\rightarrow {V\choose 2} \cup V$.  The elements of $V = V(G)$ and $E = E(G) $ are called {\em vertices and  edges of graph} $G$, respectively. If 
$\phi (e) = (v,v)$ for some $v \in V$, then $e$ is called a {\em loop of} $G$. If $\phi (a) = \phi (b)$ for some $a, b \in E$, then $a$ and $b$ are called {\em parallel edges  of} $G$.

\begin{figure}[h]
\begin{center}
\scalebox{0.3}[.3]{\includegraphics{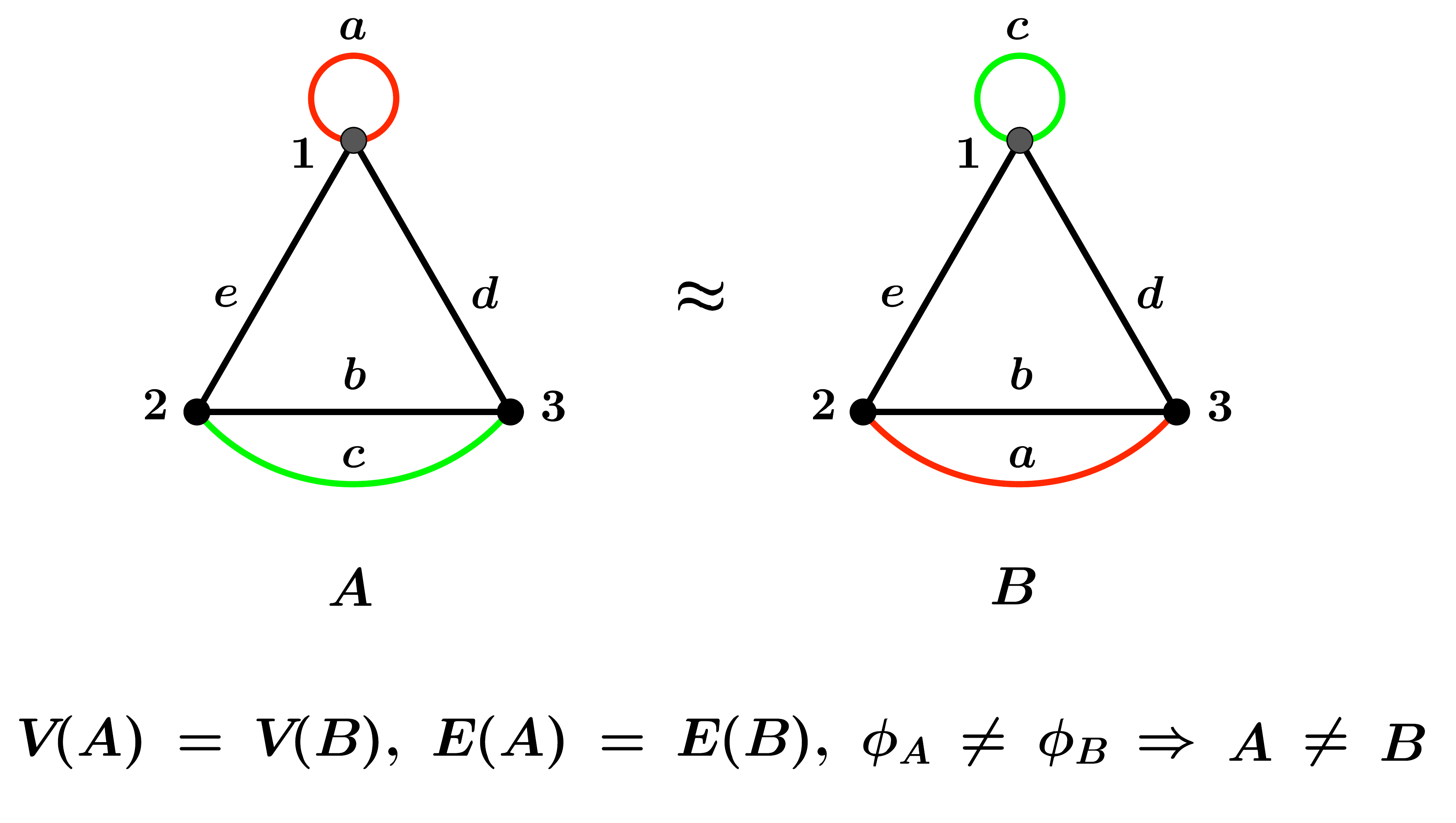}}
\end{center}
\caption
{Two different but isomorphic graphs}
\label{TsPresGdef}
\end{figure}

 \vskip 1.5ex
 Graphs  $G =( V, E ,\phi)$ and $G' = (V', E', \phi')$ are {\em equal} if 
$V = V'$, $E = E'$, and $\phi = \phi'$ (see Figure \ref{TsPresGdef}).

For graphs $G_1 = (V_1, E_1,\phi_1)$ and $G_2 = (V_2, E_2,\phi_2)$ with 
$E_1 \cap E_2 = \emptyset$, let $G$ be the graph $(V, E, \phi)$, where 
$V = V_1\cup V_2, E = E_1 \cup E_2$, and for $e \in E$, $\phi (e) = \phi_i(e)$ if 
$e \in E_i$, $i = 1, 2$. Then $G$ is called the {\em union} of $G_1$ and $G_2$, written $G_1 \cup G_2$. 
 
  \vskip 1.5ex

 We called a graph $G = (V, E, \phi)$ {\em simple} if 
 $\phi: E\rightarrow {V\choose 2} \cup V$ is an injective function, and so $G$ ha no parallel edges.

 \vskip 1.5ex

A graph $(V, E, \phi)$ is called {\em complete} if  $\phi: E\rightarrow {V\choose 2}$ is a bijection.

\vskip 1.5ex

If $\phi (e) = (u,v)$ in a graph $G$, then we say that vertices $u$ and $v$ are {\em adjacent}, edge $e$ and vertex $v$ are {\em incident in $G$}, and $u, v$ are the {\em end vertices of edge $e$ in $G$}. 
 
 \vskip 1.5ex
 A vertex in $G$ incident to no edge is called an {\em isolated vertex of} $G$. 
 
 \vskip 1.5ex
 A vertex $v$ in $G$ is called an {\em leaf of} $G$ if $v$ is incident to exactly one edge of $G$ and this edge is not a loop.
 
 \vskip 1.5ex
 
 Unless stated otherwise, we assume that a graph has no isolated vertices.
 
  \vskip 1.5ex
  
Given a graph $G = (V, E, \phi)$, an element $a \in V \cup E$, and a set 
$A \subseteq V \cup E$, we say that {\em $a$ is incident to $A$ } if $a$ is incident to an element of $A$.

 \vskip 1.5ex

\begin{figure}[h]
\begin{center}
\scalebox{0.3}[.3]{\includegraphics{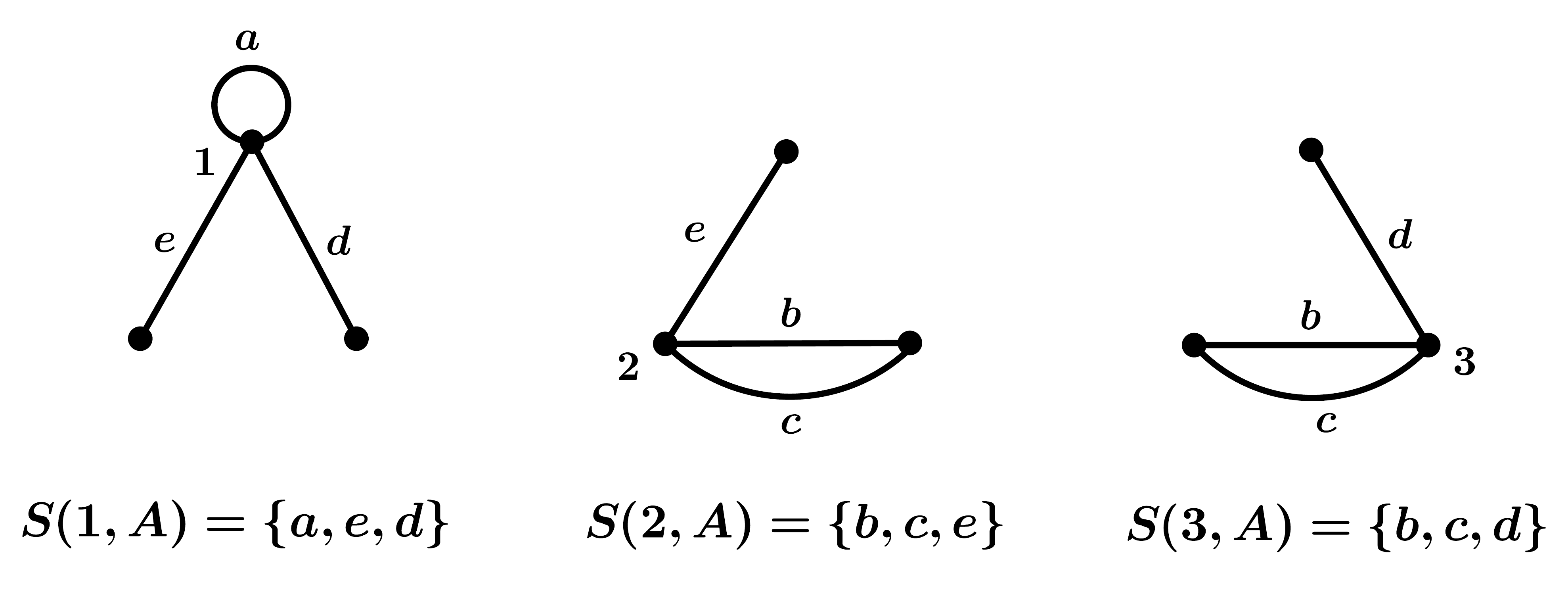}}
\end{center}
\caption
{The vertex stars of a graph}
\label{TsPresStars}
\end{figure}

 Given $G$ and $v \in V(G)$, the {\em $v$-star} (or a {\em vertex star} of $G$) is the set 
 $S(v,G)$ of edges in $G$ incident to vertex $v$. Let 
 ${\cal S}(G) = \{ S(v,G): v \in V(G)\}$ (see Figure \ref{TsPresStars}).
 
  \vskip 1.5ex
  
The {\em degree $d(v, G)$ of vertex $v$ in $G$} is the number of non-loop edges in $S(v,G)$ plus 2-times the number of loops in $S(v,G)$.
 The function  $d : V \to \N_0$ such that $d(v) = d(v, G)$ for $v \in V$ is called the {\em degree function of  $G$} and denoted by $d(G)$.
 
  \vskip 1.5ex
  
  Let  $G = (V, E,\phi)$ be a graph, $X \subset V$, and $Y \subseteq E$.
 Then $G \setminus X$ denotes the graph $G' = (V', E', \phi')$ such that
 $V' = V \setminus X$, $E'$ is obtained from $E$ by deleting all edges 
 of $G$
 incident to at least one vertex in $X$ , and $\phi': E' \to {V'\choose 2} \cup V'$ is the restriction of function $\phi $ on $E'$. We say that 
 {\em $G'$ is obtained from $G$ by deleting  vertex set $X$}.
 
 Similarly, $G \setminus Y$denotes the graph $G' = (V', E', \phi' )$ such that
 $V' = V $, $E' = E \setminus Y$, and $\phi': E' \to {V' \choose 2} \cup V'$ is the restriction of $\phi $ on $E'$. We say that {\em $G'$ is obtained from $G$ by deleting edge set  $Y$}.

  \vskip 1.5ex

A graph $G'$ is a {\em subgraph of  $G$}, written as $G' \subseteq ^s G$, 
if $G'$ can be obtained from $G$ by deleting some  edge subset $Y$ and some vertex subset $X$ of $G$. Obviously,  $ \subseteq ^s $ is a partial order on the set of graphs.

 \vskip 1.5ex

Given a graph $G = (V, E, \phi )$ and $Y \subseteq E$, let 
$G \langle Y \rangle $ denote the subgraph of $G$ obtained from $G$ by
deleting  all edges in $E \setminus Y$ and all vertices of $G$ that are incident to no edge in $Y$. We call $G \langle Y \rangle $ the {\em subgraph of $G$ induced by $Y$}.
We also say that {\em $Y$ spans $X \subseteq V$ in $G$} if 
$X = V(G \langle Y \rangle)$.
Obviously, if  $A$ is a subgraph of $G$ with no isolated vertices, 
then $G\langle  E(A) \rangle  = A$.

\vskip 1.5ex

Let ${\cal G}^{1,1}_2$ denote the set of graphs having exactly two vertices, say $x$ and $y$,  of degree $1$ and all other vertices of degree $2$. 
Then a graph $P$ is called an {\em $(x,y)$-path} (or simply, a {\em path})  if $P$ is a $\subseteq ^s$-minimal graph in  ${\cal G}^{1,1}_2$. Obviously, $x \ne y$.
We call vertices $x$ and $y$  the {\em end-vertices of path $P$} and 
put $End(P) = \{x,y\}$.

  \vskip 1.5ex

Let ${\cal G}_2$ denote the set of graphs with  all vertices having  degree two. 
Then a graph $C$ is called a {\em cycle} if $C$ is a $\subseteq ^s$-minimal graph in  ${\cal G}_2$. 
  \vskip 1.5ex

Let ${\cal G}^{1,3}_2$ denote the set of graphs having exactly one vertex, 
say $x$, of degree one,  exactly one vertex, 
say $y$, of degree three,  and all other vertices of degree two.
Then a graph $Q$ is called a {\em lollypop} if $Q$ is a $\subseteq ^s$-minimal graph in  ${\cal G}^{1,3}_2$. We call vertex $x$ the {\em end-vertex of lollypop $Q$}, put $End(Q) = \{x\}$, and also call $Q$ an {\em $x$-lollypop}.

\vskip 2.5ex
\begin{figure}[h]
\begin{center}
\scalebox{0.4}[.4]{\includegraphics{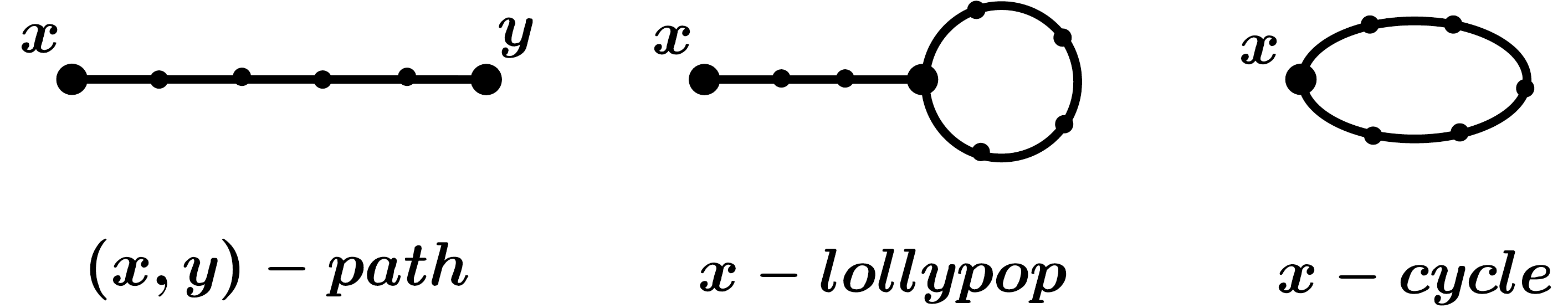}}
\end{center}
\caption
{$(x,y)$-paths, $x$-lollypops, and $x$-cycles.}
\label{Ears}
\end{figure}

  \vskip 1.5ex 
 
In a plain  language an $(x,y)$-path is a graph $P = (V, E, \phi )$ with  
$V = \{x= v_1, \ldots , v_n = y\}$ and $E = \{e_1, \ldots , e_{n-1} \}$ and
$\phi (e_i) = \{v_i, v_{i+1}\}$ for $i \in \{1, \ldots , n-1\}$.
A graph $G$ is a cycle if and only if $G$ can be obtained from an $(x,y)$-path by adding a new edge $e$ incident to $x$ and $y$.
A graph $Q$ is a lollypop if and only if $Q$ is obtained from disjoint cycle $C$
(possibly, a loop) and a path $P$  by identifying a vertex of $C$ and an end-vertex of $P$. 

\vskip 1ex

Let ${\cal C}l(G)$ denote the set of subgraphs of $G$ that are cycles. 
Let ${\cal C}t(G)$ denote the set of edge sets of elements in ${\cal C}l(G)$.

\vskip 1ex
  
  Given a graph $G$ and edge $e$ with the end vertices $x$ and $y$, we say that {\em graph $G'$ is obtained from by a subdivision of edge $e$}
  if $G' = (G \cup xPy)  \setminus e$, where $xPy$ is a path and  $V(G) \cap V(P) = \{x, y\}$. A graph $H$ is a {\em subdivision of graph $G$} if $H$ is obtained from $G$ by subdivisions of some of its edges. 

  \vskip 1ex

A {\em $k$-skein}, $k \ge 3$,  is a subdivision of a graph with two vertices and $k$ parallel edges.  A 3-skein is also called a {\em $\Theta $-graph} (see Figure \ref{BicyclesFig}).

  \vskip 1ex

A {\em dumbbell} is a graph obtained from two disjoint cycles $C$, $C'$, and a path $xPx'$  disjoint from $C$ and $C'$ by identifying vertex $x$ with a vertex in $C$ and vertex $x'$ with a vertex in $C'$ (see Figure \ref{BicyclesFig}).

  \vskip 1.5ex

A {\em butterfly} is a graph obtained from two disjoint cycles $C$, $C'$ by identifying one vertex from $C$ with one vertex in $C'$ (see Figure \ref{BicyclesFig}).

  \vskip 1ex
  
Let ${\cal G}^{3,3}_2$ denote the set of graphs having exactly two vertices
of degree three and all other vertices of degree two. Then, obviously,  
 a $\subseteq ^s$-minimal graph in ${\cal G}^{3,3}_2$ is either a 
$\Theta $-graph or a dumbbell.

  \vskip 1.5ex

Let ${\cal G}^4_2$ denote the set of graphs having exactly one vertex
of degree four and all other vertices of degree two.  Then  
 a $\subseteq ^s$-minimal graph in ${\cal G}^4_2$ is a butterfly.

\vskip 1.5ex

A graph $G$ is {\em connected}  if 
$|V(G)| \ge 2$ and 
for  every two distinct vertices $x$ and $y$ in 
$G$ there exists  an $(x,y)$-path that is  a subgraph of $G$.  

  \vskip 1.5ex
 
  A graph $G$ is {\em 2-connected}  if $|V(G)| \ge 2$, $G$ is connected and has no loops, and one of the following holds:
  
 $(a1)$ if $|V(G)| = 2$, then $|E(G)| \ge 2$ (i.e. $G$ has at least two parallel edges) and
 
 $(a2)$ if $|V(G)| \ge 3$, then $G \setminus  v$ is connected for every $v \in V(G)$.

 \vskip 1.5ex

A {\em multi-triangle} is either $K_3$ or  a graph obtained from $K_3$ by replacing some edges by parallel edges.

  \vskip 1.5ex

A graph $G$ is {\em 3-connected}  if $|V(G)| \ge 4$, $G$ has no loops 
and no parallel edges,   
and $G \setminus  v$ is 2-connected for every $v \in V(G)$.  

 \vskip 1.5ex
 
 A graph $G$ is {\em multi-3-connected} if $G$ is either 3-connected or can be obtained from a 3-connected graph by replacing some of its edges by parallel edges.

  \vskip 1.5ex

A graph $G$ is {\em k-connected} for $k \ge 4$ if $|V(G)| \ge k+1$,   and $G \setminus  v$ is $(k - 1)$-connected for every $v \in V(G)$.

 \vskip 1.5ex
 
A graph $G$ is {\em of connectivity $k$} for $k \ge 1$ if $G$ is $k$-connected but not $(k+1)$-connected. 

  \vskip 2ex

We call a graph $G$  {\em  cycle connected} if $|E(G)| \ge 2$
and for every two edges $a$ and $b$ in $G$ there exists a cycle $C$ in 
$G$ such that $a, b \in E(C)$.

  \vskip 1.5ex

A {\em component of a graph $G$} is a $\subseteq ^s$-maximal connected subgraph of $G$.
Let $Cmp (G)$ denote the set of components of $G$ and let 
$|Cmp (G)| = cmp (G)$.

  \vskip 1.5ex

A graph $G$ is called {\em cacti-graph} if G has no isolated vertices, no leaves, and  no cycle components. A connected cacti-graph is called a 
{\em cactus}. Let ${\cal G}_{\bowtie}$ denote the set of cacti-graphs and 
${\cal CG}_{\bowtie}$ denote the set of connected graphs from ${\cal G}_{\bowtie}$, and so each member of ${\cal CG}_{\bowtie}$ is a cactus
(for example, the graphs in Figures \ref{Bicycles2Fig} and \ref{TricyclesFig} and their subdivisions). Given a graph $G$, let ${\cal G}_{\bowtie} (G)$ denote the set of subgraphs of $G$ that are members of ${\cal G}_{\bowtie}$. The class ${\cal G}_{\bowtie}$ of cacti-graphs will play a special role in our further discussion. Obviously, if a graph $G$ is in ${\cal G}_{\bowtie}$, then every component of $G$ is also in ${\cal G}_{\bowtie}$. Therefore every graph that is $\subseteq ^s$-minimal in ${\cal G}_{\bowtie}$ is also $\subseteq ^s$-minimal in ${\cal CG}_{\bowtie}$.
In other words, a graph $G$ is $\subseteq ^s$-minimal in ${\cal G}_{\bowtie}$ if and only if 
 $G$ is $\subseteq ^s$-minimal in ${\cal CG}_{\bowtie}$. Below (see Claim \ref{BicyStructure}) we will show that $G$ is $\subseteq ^s$-minimal in 
 ${\cal G}_{\bowtie}$ if and only if $G$ is either a{ $\Theta $-graph or a dumbbell or a butterfly.

\vskip 2ex
\begin{figure}[h]
\begin{center}
\scalebox{0.3}[.3]{\includegraphics{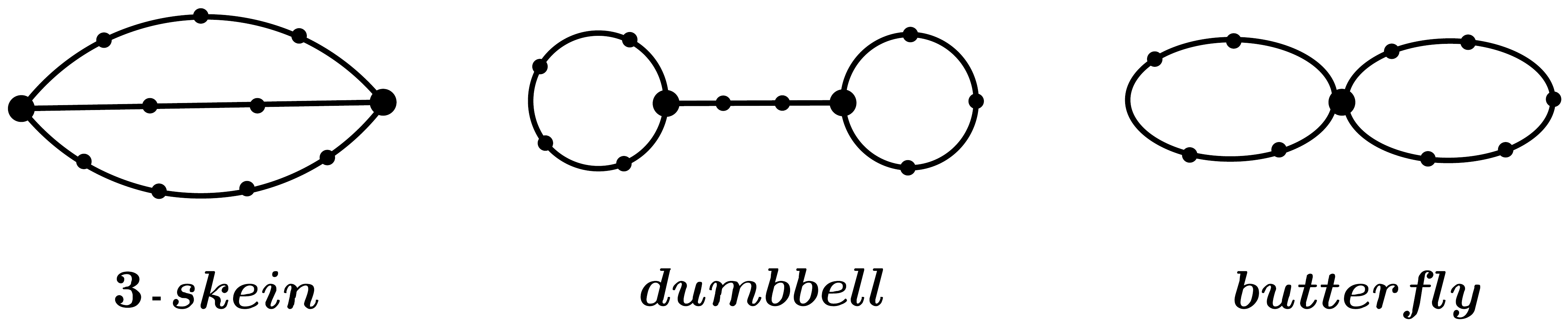}}
\end{center}
\caption
{The structure of the minimal graphs in ${\cal G}_{\bowtie}$.}
\label{BicyclesFig}
\end{figure}

A subset $X$ of vertices (of edges) in a graph $G$ is called a {\em vertex cut} 
(respectively, an  {\em edge cut}) of $G$ if $G \setminus X$ has more components than $G$.
A vertex cut (an edge cut) of $G$ consisting of one element is called a {\em cut-vertex} (respectively, a {\em cut-egde} or a {\em bridge}) of $G$.

\vskip 1.5ex

It is easy to prove the following claim.

\begin{claim}
\label{CyclConG}

Let $G$ be a graph and $G'$ a graph obtained from $G$ by deleting all isolated vertices.
Then  $G$ is cycle connected if and only if $G'$ is 2-connected. 
 
\end{claim}

A {\em tree} is a connected graph without cycles. 
A {\em wheel} $W = (V, E, \phi)$ is a simple graph
obtained from a cycle $C$ with at least three vertices by adding a new vertex $x$ (the \em center of $W$}) and the set of edges 
$\{e_v: v \in V(C) ~ and~ \phi (e_v)= \{v,x\}\}$.

 \vskip 1ex

\begin{claim}
\label{wheel}

Let  $G$ be a graph with $v(G) = v$ and $e(G) = e$.  
Suppose that $G$ is 3-connected  and has a vertex $x$ such that 
$d(x,G) > e - v$.
Then  $G$ is the wheel with center $x$, and so $d(x,G) = v - 1$ and 
$d(z,G) = 3$ for every $z \in V(G) \setminus x$.

\end{claim}

\bp Since $d(x,G) > e - v$, clearly $e(G \setminus x) < e - (e - v) = v$ and 
$v(G \setminus x) = v - 1$.
Since $G$ is 3-connected,  $G - x$ is 2-connected. Therefore 
$e(G \setminus x) = v - 1$,  $G \setminus x$ is a cycle, and vertex $x$ is adjacent to every vertex of $G \setminus x$. Thus, $G$ is the wheel with center $x$.
\ep

\vskip 1.5ex

Using the arguments similar to those in the above proof it is also easy to prove the following claim.  

\begin{claim}
\label{likewheel}

Let  $G$ be a graph.  
Suppose that $G$ is 2-connected  and has a vertex $x$ such that 
$d(x,G) > e(G) - v(G) + 1$. Then every cycle of $G$ contains $x$.

\end{claim}

From Claims \ref{wheel} and \ref{likewheel} we have: 

\begin{claim}
\label{K4wheel}

Let  $G$ be a graph with $v(G) = v$ and $e(G) = e$.  
Suppose that $G$ is 3-connected and not a complete graph on n vertices. 
Then one of the following holds: 
\\[1ex]
$(c1)$ either every vertex of $G$ has degree at most $e - v$ or 
\\[1ex]
$(c2)$ $G$ has exactly one vertex $x$ of degree $e - v + 1$ and every other vertex has degree at most $e - v$. 

\end{claim}

Graphs  $G =( V, E ,\phi)$ and $G' = (V', E', \phi')$ are {\em equal} if 
$V = V'$, 
$E = E'$, and $\phi = \phi'$.

\vskip 1.5ex

An {\em isomorphism from $G  = (V, E, \phi)$ to $G' = (V', E', \phi')$ } is a pair $(\nu, \varepsilon)$, where 
\\
$\nu: V \to V'$ and $\varepsilon: E \to E'$ are bijections such that
 $\phi(e)= \{x,y\} \Leftrightarrow \phi' (\varepsilon (e)) = \{ \nu(x),\nu(y)\}$.     
Graphs $G$ and $G'$ are {\em isomorphic} (denoted by $G \approx G'$) if there exists an isomorphism from $G$ to $G'$ (or, equivalently, an isomorphism from $G'$ to $G$).

\vskip 1.5ex

Let $G  = (V, E, \phi)$ and $G' = (V', E', \phi')$ with a bijection 
$\varepsilon: E \to E'$. Then $G$ and $G'$ are called 
{\em strongly $\varepsilon $-isomorphic } if there exists a bijection $\nu: V \to V'$  such that $(\nu, \varepsilon)$ is an isomorphism from $G$ to $G'$.

\begin{figure}[h]
\begin{center}
\scalebox{0.3}[.3]{\includegraphics{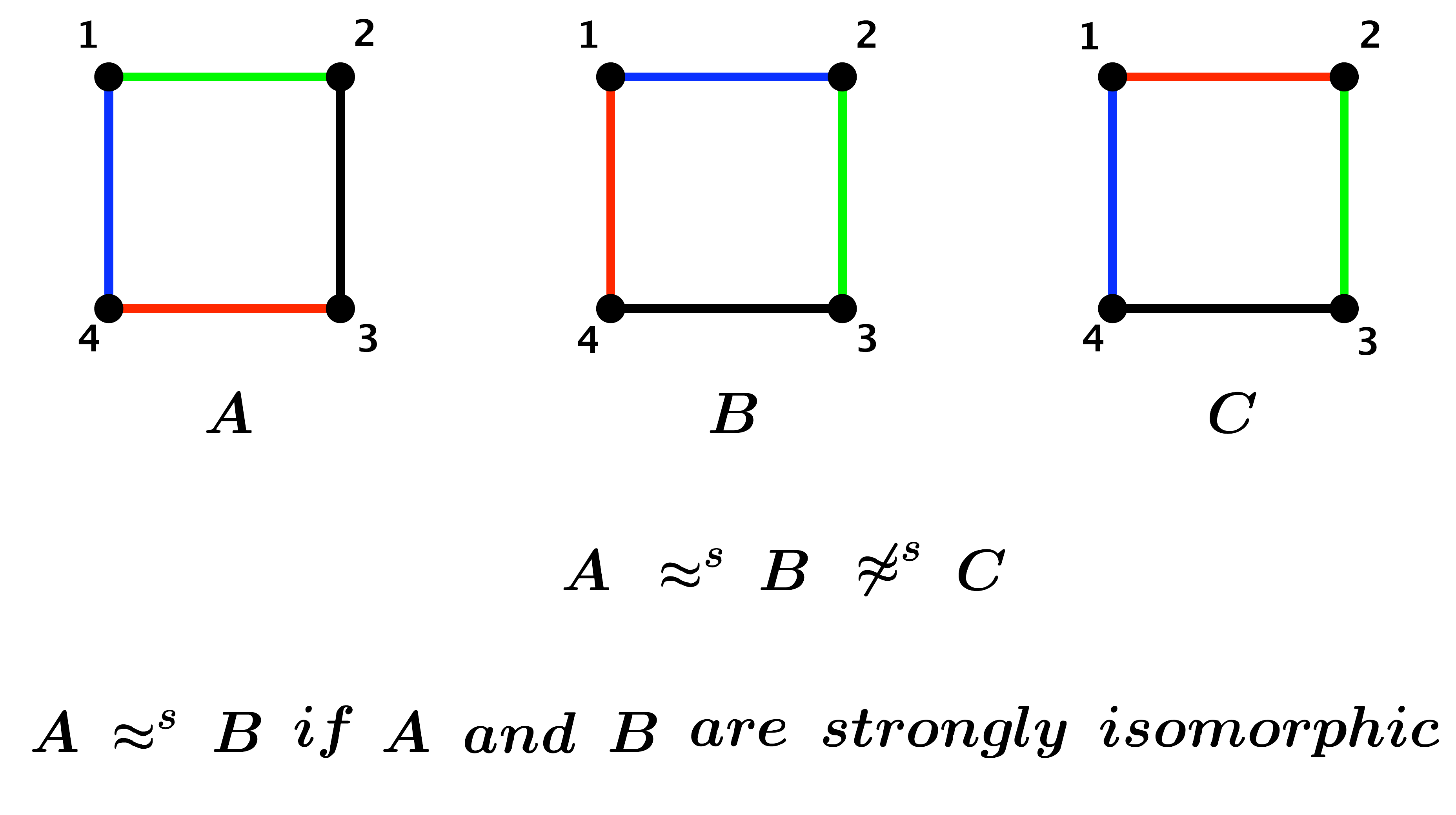}}
\end{center}
\caption
{Isomorphic and not strongly isomorphic graphs.}
\label{SquaresFig}
\end{figure}

\begin{figure}[h]
\begin{center}
\scalebox{0.3}[.3]{\includegraphics{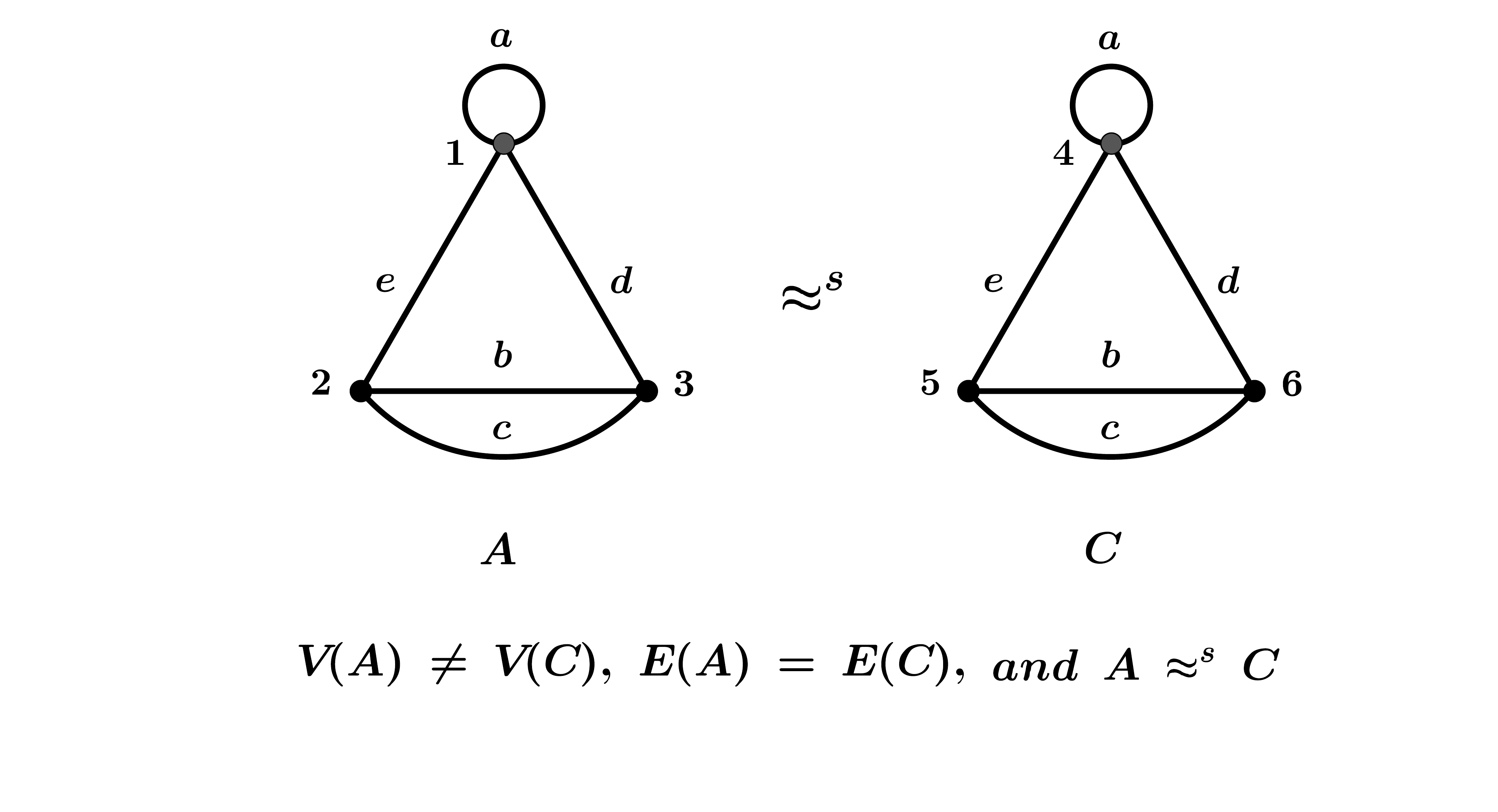}}
\end{center}
\caption
{Strongly isomorphic graphs.}
\label{TsPresStroIso}
\end{figure}

Without loss of generality, we may (and will) consider strongly isomorphic graphs  $G$ and $G'$ with $E = E'$ (i.e. when $\varepsilon: E \to E'$ is the identity function).
Namely, given two graphs $G = (V, E, \phi )$ and $G' = (V', E',\phi ')$ with 
$E = E'$,
a {\em strong isomorphism from}
$G$ to $G'$ is a bijection $\nu :V \to V'$ such that 
$\phi (e) =\{x,y\}\Leftrightarrow \phi '(e)= \{\nu(x),\nu(y)\}$. 
Graphs $G = (V, E, \phi )$ and $G' = (V', E',\phi ')$ with $E = E'$ are {\em strongly isomorphic} if
there exists a strong isomorphism from $G$ to $G'$ (see Figures \ref{SquaresFig} and \ref{TsPresStroIso}).

 \vskip 3ex
 
It is easy to prove the following claim.

\begin{claim}
\label{S(G)=S(G')}

{\em
Let $G$ and $G'$ be graphs and $E(G) = E(G')$. Then $G$ and $G'$ are strongly isomorphic if and only if ${\cal S}(G) = {\cal S}(G')$.
}
\end{claim}

\section{
Some auxiliary notions and facts on graphs}

\label{auxiliary}

\subsection{On the $\Delta$-function of a graph}
\label{}

\indent

Let ${\cal G}$ be the set of finite graphs and $G \in {\cal G}$. 
In this Section we will establish some  properties of function  
$\Delta: {\cal G} \to \Z$ such that $\Delta (G) = |E(G)| - |V(G)|$.
We will use these results on $\Delta$ in the study of so called $k$-circular matroids of  a graph that will be defined later.

\vskip 1ex
 Instead of $\Delta (G)$ we will write simply $\Delta G$. Let $X \subseteq E(G)$  
and $G\langle X \rangle $ be the subgraph of $G$ induced by $X$. Then $\Delta (G\langle  X \rangle)  = |X| - |V(G\langle X \rangle )|$.

\vskip 1.5ex

We start with the following simple observation. 

\begin{claim}
\label{DeltaCon}

Let $G$ be a connected graph. 
If $G$ has no cycle $($i.e. $G$ is a tree$)$, then $\Delta G = -1$. If $G$ has exactly one cycle, then $\Delta G = 0$. Graph $G$ has at least two cycles if and only if $\Delta G \ge 1$.

\end{claim}

It is also easy to see the following:

\begin{claim}
\label{SumDelta}

Let $G$ be a graph. Then $\Delta G = \sum \{\Delta A : A \in Cmp (G)\} $.

\end{claim}

From Claims \ref{DeltaCon} and \ref{SumDelta} we have: 

\begin{claim}
\label{NonTrivG}

Let $G$ be a graph and $t(G)$
the number of tree components of $G$. Then $G$ has a component 
with at least two cycles if and only if $
\Delta G + t(G) \ge 1$. 

\end{claim}

\begin{claim}
\label{deltaxleqdeltaa}

 Let $A$ be a connected subgraph  of $G$ and $\emptyset \neq X \subseteq E(A)$.  Then 

$\Delta G\langle  X \rangle \leq  \Delta A $. 
\end{claim}

\bp
If 
$X = E(A)$,
then our claim is obviously true. So we assume that $\emptyset \ne X \subset E(A)$. Let $S = V(A)\setminus V( G\langle  X \rangle)$. 
\\[0.7ex]
\indent
First, suppose  that $S = \emptyset $.
Since $|X| < |E(A)|$, we have:
\\
$\Delta G\langle  X \rangle = |X| - |V( G\langle  X \rangle)| < 
|E(A)| - |V(A)|= 
 \Delta A $.
Therefore our claim is true. 
\\[0.7ex]
\indent
Finally, suppose
that $S \ne \emptyset$. Since $X \ne \emptyset$ and $A$ is connected,  
there are at least $|S|$ edges of $E(A) \setminus X$ that are incident to $S$, and so $|E(A)\setminus X| \geq |S|$. Then 
\\[0.7ex]
\indent
$\Delta A =
|E(A)| - |V(A)|=
|X|+|E(A)\setminus X| - (|V( G\langle  X \rangle )| + |S|)=
\\[0.7ex]
\indent
\Delta G\langle  X \rangle + |E(A)\setminus X| - |S| \geq 
\Delta G\langle  X \rangle$. 
 \ep

\begin{claim}
\label{deltaxleqdeltah}

Let $F$ be a subgraph of  $G$  such that no component of $F$ is a tree. Then $\Delta G\langle  X \rangle \leq \Delta F$ 
for every subset $X$ of $E(F)$.

\end{claim}

\bp (uses Claim \ref{deltaxleqdeltaa})

Let $X \subseteq E(F)$ and $Cmp (F)$ be the set of components of $F$. 
\\[0.5ex] 
Let ${\cal R} = \{A \in Cmp (F): E(A) \cap X \ne \emptyset\}$.  
Then 
$\Delta G\langle  X \rangle = 
\sum_{A \in {\cal R}} \Delta G\langle E(A) \cap X\rangle $.
\\[0.5ex]  
By Claim \ref{deltaxleqdeltaa},
$\Delta G\langle E(A)\cap ~X\rangle \leq \Delta A$ for every 
$A \in {\cal R}$. 
\\[0.5ex] 
Also, since every component of $F$ is not a tree, 
$\Delta A \geq 0$ 
for every $A \in Cmp (F) \setminus {\cal R}$.
\\[0.5ex] 
Thus,~

\vskip 0.7ex 
$\Delta G\langle  X \rangle = 
\sum \{\Delta G\langle E(A)\cap ~X\rangle: {A \in{\cal R}} \} 
\leq 
 \sum \{\Delta A: {A \in {\cal R}} \} 
 \\[1.2ex]
~~~~~~~~~~~\leq 
 \sum \{ \Delta A: {A \in Cmp (F)} \}
= \Delta F.$
 \ep

\vskip 1.5ex

Now we will describe a discrete analog of the classical Intermediate Value Theorem on a continuos function.

  \vskip 1.5ex

A pair $D = (V, E)$ is called  a {\em directed graph {\em  (or simply, a {\em digraph}) if $V$ is a non-empty finite set  and
$E \subseteq \{(x,y): \{x, y\} \subseteq V~and~ x \ne y\}$. 
The elements of $V$ and $E$ are called the {\em vertices} and {\em arcs} of $D$, respectively. 
A vertex $v$ in $D$ is called {\em minimal} ({\em maximal}) if 
$(x,v) \not \in  E$ (respectively, $(v,x) \not \in  E$) for every $x \in V$.
 
  \vskip 1ex

Given $n, m \in \Z$ and $n \le m$, let $[n, m] = \{x \in \Z : n \le x \le m\}$.

\begin{Theorem} {\sc Discrete Intermediate Value Theorem}
\label{discrete-analog}

\vskip 1ex
Suppose that digraph $D = (V, E)$ and function $f: V \to \Z$ satisfy the following conditions:
\\[1ex]
$(d1)$ 
$D$ has no directed cycles, 
\\[1ex]
$(d2)$ $D$ has exactly one minimal vertex $p$
 and  exactly one maximal  vertex $b$, and  
\\[1ex]
$(d3)$
$f(p) \le f(b)$ and $|f(x) - f(y)| \le 1$ for every $(x,y) \in E$.

\vskip 1ex 
Let $R = \{f(v): v \in V\}$.
Then $ [f(p) , f(b)] \subseteq R$.

\end{Theorem}

\bp 
Suppose, to the contrary, that $ [f(p) , f(b)] \not \subseteq R$, i.e.
the exists $r \in   [f(p) , f(b)]\setminus R $. 
Obviously, $ f(p) < r <  f(b)\}$. 
Let $V' = \{v \in V: f(v) < r\}$.
Clearly, $p  \in V'$ and $b  \not \in V'$. 
Therefore $V' \ne \emptyset $. 
Let $D' =  D \setminus (V \setminus V')$, and so $V(D) = V'$. 
Since $b \not \in V'$, digraph $D'$ has a maximal vertex, say $x$.
Since $x \ne b$, clearly $x$ is not a maximal vertex in $D$. Therefore
there exists $z \in V \setminus V'$ such that $(x,z)$ is an arc in $D$.
Since $z \in V \setminus V'$, we have: $f(z) \ge r$.
If  $f(z) > r$, then $f(z) - f(x) \ge 2$. However,
$(x,z) \in E \Rightarrow  |f(x) - f(z)| \le 1$, a contradiction. Therefore
 $f(z) = r$, and so $r \in R$, a contradiction. 
\ep

\begin{claim}
\label{deletingedge}

Let $G = (V, E,\phi )$ be a graph. If $X \subseteq E$ and $e \in X$, then 

\vskip 0.7ex
$|\Delta G\langle X \rangle - \Delta G \langle X \setminus ~e \rangle | \le 1$, i.e.
$\Delta G \langle X \rangle - \Delta G \langle X \setminus ~e \rangle 
\in  \{-1, 0, 1\}$.

\end{claim}

\bp 
The graph $ G \langle X \setminus ~e \rangle $ has one fewer edge than $G\langle X \rangle$ and either zero or one or two fewer vertices. 
\ep

  \vskip 1.5ex

From Theorem \ref{discrete-analog} and Claim \ref{deletingedge} we have, 
in particular:

\begin{claim}
\label{eachi}

{\em Let $G = (V, E,\phi )$
and $\Delta G\ = \delta \ge -1$.
If at least one edge in $E$ is not a loop, then for each $i \in \{-1, \ldots ,  \delta\}$ there exists $X \subseteq E$ such that $\Delta G\langle  X \rangle = i$.}

\end{claim}

\bp (uses Theorem \ref{discrete-analog} and Claim \ref{deletingedge})

\vskip 0.6ex
Let ${\cal D} = ({\cal V}, {\cal E})$ be the digraph such that ${\cal V} = 2^E$ and
if $X, Z \in {\cal V}$, then  $(X, Z) \in {\cal E} \Leftrightarrow X \subset Z~and~ |Z \setminus X| = 1$. Obviously, digraph ${\cal D}$ has no directed cycles and has exacly one minimal vertex $p = \{ \emptyset \} $ and exactly one maximal vertex 
$b = \{E \} $. Therefore  digraph ${\cal D}$ satisfies conditions $(d1)$ and $(d2)$ of Claim \ref{discrete-analog}.
\vskip 0.6ex
Let  function $f: {\cal V} \to \Z$ be defined as follows:
$f(X) = \Delta G\langle  X \rangle$ for 
$X \subseteq E$ and $X \ne \emptyset $, and so $f(E) = \delta$. 
Put $f(\emptyset ) = -1$. 
\vskip 0.6ex
Obviously,  for $e \in E$ we have: 
$f(\{e\}) = -1$ if $e$ is not a loop and $f(\{e\}) = 0$ if $e$ is a loop.
Hence in both cases 
$|f(\{e\}) - f( \emptyset )| \le  1$.
Also if $X \subset E $ and $|X| \ge 2$, then by Claim \ref{deletingedge}, 
$|f(X) - f(X \setminus e)| \le 1$. Thus, function $f: {\cal V} \to \Z$ satisfies condition $(d3)$ in Theorem \ref{discrete-analog}. Therefore our claim follows from Theorem \ref{discrete-analog}.
 \ep

\subsection{Splitter theorems for cacti-graphs}
\label{splitter}

\indent

First we will remind the following known theorem 
which is very easy to prove.

\begin{Theorem} {\sc Splitter Theorem for 2-connected graphs}
\label{EarDcmp2con}

Let  $G$ be a $2$-connected graph and $G_0$ a proper 2-connected  subgraph of $G$.
Then  there exists a sequence
$(P_1, \ldots , P_r)$ of paths in $G$ and a sequence $(G_0, G_1, \ldots , G_r)$  of 2-connected subgraphs of $G$ such that  $G_r = G$, and 
$G_i = G_{i-1}\cup P_i$, where $V( G_{i-1}) \cap V(P_i) = End(P_i)$ for every $i \in \{1, \ldots , r\}$, and so  $G_{i-1} \subset ^s G_i$ and
$\Delta G_i - \Delta G_{i-1} = 1$.

\end{Theorem}

As in Section \ref{OnGraphs},  ${\cal G}_{\bowtie}$ (${\cal CG}_{\bowtie}$) is the set of graphs (respectively, connected graphs) with no isolated vertices, no leaves, and  no cycle components, i.e. the set of cacti-graphs, (respectively, connected cacti-graphs). Notice that if graph $G$ is 2-connected and it is not a cycle, then $G \in {\cal CG}_{\bowtie}$.
\vskip 1.5ex
It turns out that  similar  Splitter Theorems are also true for the graph classes  
${\cal CG}_{\bowtie}$ and ${\cal G}_{\bowtie}$. These theorems will be important for our further  discussion.

\vskip 1.5ex

\begin{Theorem} {\sc Splitter theorem for a cactus}
\label{EarDcmpConGinfty}

Let  $G$ and $G_0$ be connected graphs. Suppose that  $G \in {\cal CG}_{\bowtie}$ 
and $G_0$ is a  subgraph of $G$ such that either 
$G_0\in {\cal CG}_{\bowtie}$ or $G_0$ is a cycle.
Then  there exist  sequences
$(P_1, \ldots , P_r)$ and 
$(G_0, G_1, \ldots , G_r)$  of subgraphs of  $G$ such that $G_r = G$, 
every $G_i \in {\cal CG}_{\bowtie}$, every $P_i$ is either a path or a  cycle or a  lollipop {\em(see Figure \ref{Ears})}, and for every 
$i \in \{1, \ldots , r\}$, we have: $G_i = G_{i-1}\cup P_i$, where 
\\[1ex]
${\bf (p)}$ $V( G_{i-1}) \cap V(P_i) = End(P_i)$ if $P_i$ is a  a path,
\\[1ex]
${\bf (l)}$ $V( G_{i-1}) \cap V(P_i) = End(P_i)$ if $P_i$ is a lollipop, and
\\[1ex] 
${\bf (c)}$ $|V( G_{i-1}) \cap V(P_i)| = 1$ if $P_i$ is a cycle, 
\\[1ex]  
and so  
$G_{i-1} \subset ^s G_i$ and 
$\Delta G_i - \Delta G_{i-1} = 1$.

\end{Theorem}

We call the sequence 
$(G_0 \subset ^s  G_1 \subset ^s  \ldots \subset ^s G_r = G)$ 
in the above  theorem
an {\em ear-assembly of a  graph $G$ in ${\cal G}_{\bowtie}$ starting from $G_0$} and each $P_i$ an {\em ear of this ear-assembly}.
\\[1.5ex]
{\bf Proof of Theorem \ref{EarDcmpConGinfty}.}
We will first consider the case when  $G_0 \in {\cal CG}_{\bowtie}$.
The case when $G_0$ is a cycle will follow easily from the first case.
\\[1ex]
${\bf (p1)}$ Suppose that $G_0 \in {\cal CG}_{\bowtie}$.
Suppose, to the contrary, that our claim is not true. Our theorem is vacuously true if $G = G_0$.
Suppose that $G_k$ is a maximal subgraph of $G$ such that our claim is true 
for  $G_k$, and so there exist   sequences $(P_1, \ldots , P_k)$ and 
$(G_0,   \ldots , G_k)$  satisfying all the conditions of our theorem.
Then $G_k \in {\cal CG}_{\bowtie}$ and $G_k$ is a proper subgraph of $G$
$($i.e. $G_k \subset ^s G$$)$.
We will get a contradiction by showing that  there exist subgraphs 
 $P_{k+1}$ and $G_{k+1}$ of $G$ such that $G_k \cup P_{k+1} = G_{k+1}$,  
$G_{k+1} \in  {\cal G}_{\bowtie}$, and pair $(G_k, P_{k+1})$  satisfies  one of conditions ${\bf (p)}$,${\bf (l)}$, ${\bf (c)}$.

Since $G$ is connected and $G_k$ is a proper subgraph of $G$, there is a path $xPy$ in $G$ such that $x \in V(G_k)$ and $V(P)\cap V(G_k) \subseteq \{x,y\}$. We can assume that $xPy$ is a maximal path with the above mentioned properties.
\vskip 0.6ex
First, suppose that $V(P)\cap V(G_k) = \{x,y\}$. 
Put $xPy = P_{k+1}$ and $G_{k+1} = G_k \cup  P_{k+1}$. 
Then $G_{k+1} \in {\cal CG}_{\bowtie}$, $G_{k+1} \subseteq ^s G$, 
and pair $(G_k, P_{k+1})$  satisfies  condition ${\bf (p)}$.
\vskip 0.6ex
Finally, suppose that $V(P)\cap V(G_k) = x$. Since $G \in {\cal CG}_{\bowtie}$, vertex 
$y$ is not a leaf in $G$.  Hence there exists edge $e = yz$ in
 $G \setminus (G_k \cup P)$. Since  $xPy$ is a maximal path with the above mentioned properties, vertex $z$ is in $P$. 

If $z \ne x$, then $P \cup e$ is a lollipop. Put $ P_{k+1} = P \cup e$ and $G_{k+1} = G_k \cup  P_{k+1}$. 
Then $G_{k+1} \in {\cal CG}_{\bowtie}$, $G_{k+1} \subseteq ^s G$, 
and pair $(G_k, P_{k+1})$  satisfies  condition ${\bf (l)}$.

If $z = x$, then $P \cup e$ is a cycle. Put $ P_{k+1} = P \cup e$ and 
$G_{k+1} = G_k \cup  P_{k+1}$. 
Then $G_{k+1} \in {\cal CG}_{\bowtie}$, $G_{k+1} \subseteq ^s G$, 
and pair $(G_k, P_{k+1})$  satisfies  condition ${\bf (c)}$.
\\[1ex]
${\bf (p2)}$
Now suppose that $G_0$ is a cycle. The arguments similar to those in 
${\bf (p1)}$ show that there exist   subgraphs  $P_1$ and 
$G_1$  of $G$ such that $G_0 \cup P_1 = G_1$,  
$G_1 \in  {\cal CG}_{\bowtie}$, and pair $(G_0, P_1)$  satisfies  one of conditions ${\bf (p)}$,${\bf (l)}$, ${\bf (c)}$.
Now put in ${\bf (p1)}$ $G_0: = G_1$. Then by  ${\bf (p1)}$, our theorem is true  for $G_0: = G_1$, i.e. there exists ear assembly $(G_1 \subset ^s  G_2 \subset ^s  \ldots \subset ^s G_r = G)$. Then 
$(G_0 \subset ^s  G_1 \subset ^s  \ldots \subset ^s G_r = G)$ is an ear assembly in ${\cal CG}_{\bowtie}$ of $G$ starting from $G_0$.   
\ep

\vskip 1.5ex

In order to formulate and prove a Splitter Theorem for ${\cal G}_{\bowtie}$,
 we need some preliminaries.

\begin{definition}
\label{bicycle}
{\em
A graph $D$ is a {\em bicycle} if $D$ is a connected graph with no leaves and with $\Delta D = 1$. 
}
\end{definition}

Now we can describe the structure of graph-bicycles.

\begin{claim}
\label{BicyStructure}

A graph $D$ is a bicycle if and only if $D$ is either a $\Theta$-graph or a dumbbell or a butterfly {\em(see Figure \ref{BicyclesFig})}. 
 
\end{claim}

\bp (uses Theorem \ref{EarDcmpConGinfty})

Obviously, a $\Theta$-graph, a dumbbell, and a butterfly is a bicycle.
It remains to prove that every bicycle $D$ is either a $\Theta$-graph or dumbbell or a butterfly. 
Since $D$ is a bicycle, $D$ is connected, has no leaves, and  
$\Delta D = 1$. Therefore  $D \in   {\cal CG}_{\bowtie}$ and $\Delta D = 1$, and so $D$ has a cycle $C$ and $\Delta C = 0$.  
Therefore our claim follows from Theorem \ref{EarDcmpConGinfty}, when $G_0$ is a cycle, $G_1 = D$, and $r = 1$.
\ep

\vskip 1.5ex

It turns out that  every graph $G$ in ${\cal G}_{\bowtie}$  has bicycles and, moreover, graph $G$ has some special bicycles.

\begin{claim}
\label{AhasBicy}
 
Suppose that $G$ is a   graph in ${\cal CG}_{\bowtie}$. Then for every two edges $a$, $b$ in $G$ there exists a bicycle in $G$ containing $a$ and $b$. 
 
\end{claim}

\bp 
By Claim \ref{BicyStructure}, it is sufficient to prove that $G$ has either a $\Theta$-graph or a dumbbell or a butterfly containing $a$ and $b$. 
\\[1ex]
${\bf (p1)}$
Suppose that $G$ has two cycles $C$, $C'$ such that $a, b \in E(C \cup C')$. We can assume that $\{a, b\} \cap E(C) \ne \emptyset$. 

Suppose first that $|V(C) \cap V(C')| \ge 2$. Then in $C'$ there exists a path $P$ such that $V(C) \cap V(P) = End(P)$ and $D = C \cup P$ contains both $a$ and $b$. Then $D$ is a $\Theta$-graph in $G$ containing $a$ and $b$. 

Now suppose that $|V(C) \cap V(C')| \le 1$. Since $G$ is connected, there exists a path $xPx'$ in $G$ (possibly, $x = x'$) such that $V(C) \cap V(P) = x$ and $V(C') \cap V(P) = x'$. If $|V(P)| = 1$, then $D = C \cup C'$ is a butterfly in $G$ containing $a$ and $b$. If $|V(P)| \ge 2$, then $D= C \cup P \cup C'$ is a dumbbell in $G$ containing $a$ and $b$.
\\[1ex]
${\bf (p2)}$
Next suppose that $G$ has a cycle $C$ containing exactly one of $a, b$, say 
$a$, and $b$ belongs to no cycle in $G$. 
Since $G$ is connected, there exists a path $xPx'$ such that $b \in E(P)$ and $V(C) \cap V(P) = x$. Since 
$G \in {\cal G}_{\bowtie}$, clearly $G$ has no leaves, and therefore $P$ belongs to a lollypop $L$ in $G$ such that $V(C) \cap V(L) = End (L)$. 
Since $b$ belongs to no cycle in $G$, graph $C \cup L$ is a dumbbell containing both $a$ and $b$.   
\\[1ex]
${\bf (p3)}$ 
Finally, suppose that neither $a$ nor $b$ belongs to a cycle in $G$. Since $G$ is connected, there exists a path $P$ in $G$ containing $a$ and $b$. Since 
$G \in {\cal G}_{\bowtie}$, clearly $G$ has no leaves, and therefore $P$ belongs to two lollypops $L$ and $L'$ in $G$. Then $L \cup P \cup L'$ is a dumbbell containing $a$ and $b$.      
\ep

\vskip 1.5ex

Now we will  formulate and prove a Splitter Theorem for graph class 
${\cal G}_{\bowtie}$.

\vskip 1ex

\begin{Theorem} {\sc Splitter theorem for a cacti-graph}
\label{EarDcmpGinfty}

Let  $G, G_0 \in {\cal G}_{\bowtie}$. 
Suppose that $G_0$ is a 
subgraph of $G$.
Then  either $G = G_0$ or there exist  sequences
$(P_1, \ldots , P_r)$ and 
$(G_0, G_1, \ldots , G_r)$  of subgraphs of  $G$ such that $G_r = G$, 
every $G_i \in {\cal G}_{\bowtie}$, every $P_i$ is either a path or a cycle or a lollipop {\em(see Figure \ref{Ears})} or a bicycle, and for every 
$i \in \{1, \ldots , r\}$, we have: $G_i = G_{i-1}\cup P_i$, where 
\\[1ex]
${\bf (p)}$ $V( G_{i-1}) \cap V(P_i) = End(P_i)$ if $P_i$ is a  a path,
\\[1ex]
${\bf (l)}$ $V( G_{i-1}) \cap V(P_i) = End(P_i)$ if $P_i$ is a lollipop,
\\[1ex] 
${\bf (c)}$ $|V( G_{i-1}) \cap V(P_i)| = 1$ if $P_i$ is a cycle, and 
\\[1ex]  
${\bf (b)}$ $V( G_{i-1}) \cap V(P_i) = \emptyset$ if $P_i$ is a bicycle, 
\\[1ex]  
and so  
$G_{i-1} \subset ^s G_i$ and 
$\Delta G_i - \Delta G_{i-1} = 1$.

\end{Theorem}

\bp (uses Theorem \ref{EarDcmpConGinfty})

Suppose, to the contrary, that our claim is not true. Our theorem is vacuously true if $G = G_0$.
Suppose that $G_k$ is a maximal subgraph of $G$ such that our claim is true 
for  $G_k$, and so there exist   sequences $(P_1, \ldots , P_k)$ and 
$(G_0,   \ldots , G_k)$  satisfying all the conditions of our theorem.
% and so $G_k \in  {\cal G}_{\bowtie}$. 
Then $G_k \in {\cal G}_{\bowtie}$ and $G_k$ is a proper subgraph of $G$
$($i.e. $G_k \subset ^s G$$)$.
We will get a contradiction by showing that  there exist subgraphs 
 $P_{k+1}$ and $G_{k+1}$ of $G$ such that $G_k \cup P_{k+1} = G_{k+1}$,  
$G_{k+1} \in  {\cal G}_{\bowtie}$, and pair $(G_k, P_{k+1})$  satisfies  one of conditions ${\bf (p)}$,${\bf (l)}$, ${\bf (c)}$, ${\bf (b)}$.
%
%Obviously, if graph $F$ is  in $ {\cal G}_{\bowtie}$, then every component of $F$ is in  $ {\cal G}_{\bowtie}$. 
\\[1.5ex]
${\bf (p1)}$
Suppose that $G$ has a component $A$ disjoint from $G_k$.
Since $G \in  {\cal G}_{\bowtie}$, also $A \in  {\cal G}_{\bowtie}$.
By Claim \ref{AhasBicy},  $A$ has a subgraph-bicycle $B$.
Put $P_{k+1} = B$ and $G_{k+1} = G_k \cup  P_{k+1}$. 
Then $G_{k+1} \in {\cal G}_{\bowtie}$, $G_{k+1} \subseteq ^s G$, 
and pair $(G_k, P_{k+1})$  satisfies  condition ${\bf (b)}$.
\\[1.5ex]
${\bf (p2)}$
Now suppose that $G$ has no  component disjoint from $G_k$.
Since $G_k$ is a proper subgraph of $G$,  
there are components $A$ of $G$ and $A_k$ of $G_k$ such that $A_k$ is a proper subgraph of $A$.
Since $G, G_k \in  {\cal G}_{\bowtie}$, also $A,  A_k \in  {\cal G}_{\bowtie}$. 
Therefore $A,  A_k \in  {\cal CG}_{\bowtie}$. 
Now our claim follows from 
Theorem \ref{EarDcmpConGinfty} applied to $A$ and $A_k$.
\ep

\vskip 1.5ex

We call the sequence 
$(G_0 \subset ^s  G_1 \subset ^s  \ldots \subset ^s G_r = G)$ 
in the above  theorem
an {\em ear-assembly of a  graph $G$ in ${\cal G}_{\bowtie}$ starting from $G_0$} and each $P_i$ an {\em ear} of this ear-assembly.

\vskip 1.5ex

From  Splitter Theorem \ref{EarDcmpGinfty} and   Claim \ref{AhasBicy} we also have the following useful specification  of this  theorem. 

\begin{Theorem}  
\label{a,bEarDcmpGinfty}

Let  $G \in {\cal G}_{\bowtie}$ and $a$, $b$ edges in $G$. 

If $a$ and $b$ are in the same component of $G$, then 
either $G = G_0$ or 
there exists an ear-assembly $(G_0\subset ^s G_1 \subset ^s \ldots \subset ^s G_r = G)$ of graph $G$ such that 
 $a$ and $b$ are  edges of   $G_i$, $G_i \in {\cal G}_{\bowtie}$,  and  $\Delta G_i = i +1$ for every $i \in \{0, \ldots , r\}$, and so, in particular, $G_0$ is a bicycle.

If $a$ and $b$ are in different  components of $G$, then
either $G = G_1$ or 
there exists an ear-assembly $(G_1\subset ^s G_2 \subset ^s \ldots \subset ^s G_r = G)$ of graph $G$ such that 
 $a$ and $b$ are  edges of   $G_i$, $G_i \in {\cal G}_{\bowtie}$,  and  $\Delta G_i = i +1$ for every $i \in \{1, \ldots , r\}$.

\end{Theorem}

\vskip 6ex

\begin{corollary} 
\label{a,bInFsbgrG}

Let  $G \in {\cal G}_{\bowtie}$ and $a, b \in E(G)$. 

If $a$ and $b$ are in the same component of $G$, then
for every $k \in \{1, \ldots  , \Delta G\}$  there exists a subgraph $F$ of $G$ such that 
$a, b \in E(F)$,
$F \in  {\cal G}_{\bowtie}$, and  $\Delta F = k$.

If $a$ and $b$ are in different  components of $G$, then
for every $k \in \{2, \ldots  , \Delta G\}$  there exists a subgraph $F$ of $G$ such that 
$a, b \in E(F)$,
$F \in  {\cal G}_{\bowtie}$, and  $\Delta F = k$.

\end{corollary}

Here is another characterization of  bicycles.
\begin{claim}
\label{BicyMin}

{\em A graph $D$ is a bicycle if and only if $D$ is a $\subseteq^s$-minimal graph in ${\cal G}_{\bowtie}$.
}
\end{claim}

\bp (uses Theorem \ref{EarDcmpGinfty} and Claim \ref{AhasBicy}) 
\\[1ex]
${\bf (p1)}$ 
Suppose that $D$ is a bicycle. By Claim \ref{BicyStructure}, $D$ is either a 
$\Theta$-graph or a dumbbell or a butterfly. It is easy to check that every such graph is in ${\cal G}_{\bowtie}$. Therefore $D \in {\cal G}_{\bowtie}$.

We prove that $D$ is a $\subseteq^s$-minimal graph in ${\cal G}_{\bowtie}$. Suppose not. Then there exists graph $D_0$ in  ${\cal G}_{\bowtie}$ such that $D_0$ is a proper subgraph of  $D$. By Theorem \ref{EarDcmpGinfty}, there exists  an ear assembly 
$(D_0 \subset ^s D_1  \subset ^s \ldots  \subset ^s D_r = D)$, where  
$\Delta D_i - \Delta D_{i-1} = 1$ for every $i \in \{1, \ldots , r\}$.
Therefore $\Delta D_0 < \Delta D$. Since $ \Delta D = 1$, we have: 
$\Delta D_0 < 1$. Hence $D_0$ has a component with at most one cycle implying that 
$D_0 \not \in {\cal G}_{\bowtie}$, a contradiction.
Therefore $D$ is a $\subseteq^s$-minimal graph in 
${\cal G}_{\bowtie}$. 
\\[1ex]
${\bf (p2)}$
Now suppose that $D$ is a $\subseteq^s$-minimal graph in 
${\cal G}_{\bowtie}$. Let $A$ be a component of $D$. Since $D \in {\cal G}_{\bowtie}$, also $A \in {\cal G}_{\bowtie}$. By Claim \ref{AhasBicy}, graph $A$ contains a bicycle $Q$ which is a subgraph of $D$. Since $D$ is a $\subseteq^s$-minimal graph in 
${\cal G}_{\bowtie}$, we have: $D = Q$, and so $D$ is a bicycle.    
\ep

\subsection{On the core and the kernel of a graph}

\indent

In this Section we will give a definition and 
establish some 
properties of the {\em core}  and the {\em kernel} of  graphs. 

\vskip 1.5ex
For a graph $A$, let ${\cal R}_A = \{ X \subseteq E(A) : X \ne \emptyset~ and~
\Delta A \langle X \rangle  = \Delta A \}$.
\begin{claim}    
\label{coreconnected}

Let $A$ be a connected graph with at least one cycle. 
Then 
${\cal R}_A $
has 
\\
the $\subseteq $-minimum element 
{\em (denoted by)} ${\cal M}in ({\cal R}_A)$.

\end{claim} 

\bp   (uses Claim~\ref{deltaxleqdeltaa} and  \ref{deltaxleqdeltah})

Let $A_{0} = A$ and define recursively $A_s$, for $s \geq 1$, to be the graph obtained from $A_{s-1}$ by deleting all leaves in $A_{s-1}$. 
Consider $r \geq 0$ such that $A_r$ has no leaves. 
Obviously, $A_r$ is a connected graph, $E(A_r) \ne \emptyset$  and $\Delta A_r = \Delta A$. 
We claim that $E(A_r)$ is the $\subseteq $-minimum element 
 in ${\cal R}_A $.
\indent
Suppose, to the contrary, that $E(A_r)$ is not the $\subseteq $-minimum element in ${\cal R}_A$, i.e.
there exists $X \in {\cal R} _A$ such that 
$E(A_r) \setminus  X \ne  \emptyset $. Since  $X \in {\cal R} _A$, we have: 
$ \emptyset \ne X \subseteq E(A)$  and $\Delta A \langle X \rangle =
 \Delta A$. 
Since $E(A_r) \setminus  X \ne  \emptyset $, there exists $e \in E(A_r) \setminus X$.
If the graph $A \setminus e$ is connected, then by Claim \ref{deltaxleqdeltaa},
$\Delta A \langle X \rangle \leq 
\Delta (A \setminus e) <  \Delta A$.
If $A \setminus e$ is not connected, then $A \setminus e$ has two components. Since $e \in E(A_r)$, both components of $A \setminus e$ have cycles. Therefore, by Claim \ref{deltaxleqdeltah},
$\Delta A \langle X \rangle \leq 
\Delta (A \setminus e) <  \Delta A$. 
 In both cases, $\Delta A \langle X \rangle < \Delta A$, a contradiction.   
\ep

\vskip 1.5ex

\begin{remark}
\label{T,NoMin}
{\em
If a graph $A$ is a tree with at least two edges, then ${\cal R}_A $ has no $\subseteq $-minimum element.    
}
\end{remark}

\begin{definition} {\sc The kernel of a graph}
\label{ConGraphCore}

\vskip 0.7ex
{\em
If $A$ is a connected graph with at least one cycle, 
then let $\lfloor A \rfloor$ denotes 
the subgraph of $A$ induced by the edge subset ${\cal M}in ({\cal R}_A)$ defined in Claim \ref{coreconnected}, and so 
${\cal M}in ({\cal R}_A) = E \lfloor A \rfloor $. 
\\
If $F$ is a non-connected graph, then put $  \lfloor F \rfloor = \cup \{ \lfloor A \rfloor : A \in Cmp(G)$ and  $\Delta A \ge 0 \}$. 
\vskip 0.7ex
Graph  $\lfloor G \rfloor$ is called the {\em kernel of graph $G$}. If $F$ is a forest, then the kernel of $F$ is not defined.  
}
\end{definition}

\begin{definition} {\sc The core of a graph}
\label{NonConGraphCore}
 
 \vskip 0.7ex
{\em  
Let $G$ be a graph. Suppose that $G$ has a component with at least two cycles. Put $ [ G]  = \cup \{ \lfloor A \rfloor : A \in Cmp(G) $ and 
$ \Delta A \ge 1 \}$. 
Graph $ [G] $ is called {\em the core of graph} $G$. If every component of $G$ has at most one cycle, then the core of $G$ is not defined.
}
\end{definition}

For a graph $A$, let ${\cal Q}_A = \{X \subseteq E(A) : A \langle X \rangle~ has~ no~ leaves \}$. 
\begin{claim}
\label{minmax}

Let $A$ be a connected graph with at least one cycle. Then
${\cal Q}_A$
has 
\\
the $\subseteq $-maximum element ${\cal M}ax ({\cal Q}_A)$. Moreover, ${\cal M}ax ({\cal Q}_A) = {\cal M}in ({\cal R} _A)$ and  
$C$ is a subgraph of $\lfloor A \rfloor$ for every cycle $C$ in $A$. 

\end{claim}   

\bp (uses Claim \ref{coreconnected})

Let $A_r$ be the graph defined in the proof of Claim \ref{coreconnected}. In that proof we have shown  that $ \lfloor A \rfloor = A_r $. 
It is clear that $A_r$ has no leaves. Also, any subgraph of $A$ containing $E(A_r)$ and at least one edge in $E(A) \setminus E(A_r)$ contains a leaf. Therefore, $E(A_r)$ is a $\subseteq $-maximal element in 
${\cal Q}_A$.
It is also clear that any set that is an element  in ${\cal Q}_A$
is a subset of  $E(A_r)$. 
Hence $E(A_r)$ is the $\subseteq $-maximum element in ${\cal Q}_A$.  
 \ep 

\vskip 1.5ex 

 Let, as above,
 \\[0.7ex]
  ${\cal R}_G = \{ X \subseteq E(G) : X \ne \emptyset~ and~ \Delta G \langle X \rangle  = \Delta G \}$ and 
 \\[0.7ex]
${\cal Q}_G = \{X \subseteq E(G) : G \langle X \rangle~ has~ no~ leaves \}$. 

\begin{claim}
\label{GMinMax}

Let $G$ be a non-connected graph such that every component of 
$G$ has at least two cycles. Then 
$$
{\cal M}in ({\cal R}_G)  = 
\cup \{ {\cal M}in ({\cal R}_A): A \in Cmp (G)\} =
$$
$$
\cup \{{\cal M}ax ({\cal Q}_A): A \in Cmp (G)\}  =
{\cal M}ax ({\cal Q}_G)= E [ G ].
$$

\end{claim}   

\bp (uses Claims \ref{SumDelta} and \ref{minmax})

The first equality follows easily from Claim \ref{SumDelta}. The second equality follows from Claim \ref{minmax}. The last equality follows from the fact that a graph has no leaf if and only if every component of the graph has no leaf. 
\ep

\vskip 1.5ex 

We recall that ${\cal G}_{\bowtie}$ is the set of graphs with no isolated vertices, no leaves, and  no cycle components. Given a graph $G$, let ${\cal G}_{\bowtie} (G)$ denote the set of subgraphs of $G$ that are members of ${\cal G}_{\bowtie}$.

\vskip 1.5ex
From Claims \ref{minmax} and \ref{GMinMax} we have:

\begin{Theorem} {\sc Description of the core of a graph} 
\label{CoreInGinfty}

\vskip 0.6ex 
Let $G$ be a graph having a component with at least two cycles. The following are equivalent:
\\[1ex]
$(a1)$ $F$ is the core of $G$, 
\\[1ex]
$(a2)$ $F$ a unique $\subseteq^s$-maximum element in 
${\cal G}_{\bowtie} (G)$, and 
\\[1ex]
$(a3)$ if $G'$ is the graph obtained from $G$ by removing all tree components, then $F$ is the only element in ${\cal G}_{\bowtie} (G)$ $($or, equivalently, in 
${\cal G}_{\bowtie} (G')$$)$ such that $\Delta F = \Delta G' $. 

\end{Theorem}

\begin{remark}
\label{G,NoMin}
{\em
Let $G$ be a graph. Then
\\[1ex]
$(c1)$ if $G$ has a tree component, then 
${\cal R}_G$ has no $\subseteq $-minimum element,
\\[1ex]
$(c2)$ if $G$ is not connected, $G$ has no tree component, and every component of $G$ has exactly one cycle, then again ${\cal R}_G$ has no $\subseteq $-minimum element, and  
\\[1ex]
$(c3)$ if $G$ has no tree component but $G$ has a component with exactly one cycle and also a component with at least two cycles, then 
 ${\cal M}in ({\cal R}_G)$ exists and 
 ${\cal M}in ({\cal R}_G) \subset {\cal M}ax ({\cal Q}_G) $. 
\vskip 1ex
 Thus, Claim \ref{GMinMax} provides a necessary and sufficient conditions for 
${\cal M}in ({\cal R}_G) = {\cal M}ax ({\cal Q}_G) $. 
}
\end{remark}

\begin{claim}
\label{A-e}

Let $A$ be a connected graph with at least one cycle and 
$e \in E(A)$.
Then
\\[1ex]
$(c0)$ $e \in E \lfloor A \rfloor $ if and only if both end vertices of $e$ belong to $ \lfloor A \rfloor $, 
\\[1ex]
$(c1)$  if $e \notin E \lfloor A \rfloor $, then $A \setminus e$ has  two components and exactly one of them is a tree and the other component contains $\lfloor A \rfloor $,
\\[1ex]
$(c2)$  if $A$ has one cycle and $e \in E \lfloor A \rfloor $, 
then $A \setminus e$  is a tree,
\\[1ex]
$(c3)$ if $A$ has at least two cycles and  $e \in E \lfloor A \rfloor $, then  
every component of $A \setminus e$ contains a cycle, and
\\[1ex]
$(c4)$ if $A \setminus e$ has two components and 
$v \in V \lfloor A \rfloor $, then the component of $A \setminus e$ containing $v$ is not a tree.

\end{claim}   

\bp (uses Definition \ref{ConGraphCore} and Claims \ref{deltaxleqdeltaa} and \ref{coreconnected}) 
\\[1.5ex]
${\bf (p0)}$ 
We prove $(c0)$.
Clearly, if $e \in E \lfloor A \rfloor $, then both end vertices of $e$ belong to 
$ \lfloor A \rfloor $. Now suppose that both end vertices of $e$ belong to 
$ \lfloor A \rfloor $ but $e \notin E \lfloor A \rfloor $. Then 
$ \Delta (\lfloor A \rfloor \cup e) = \Delta \lfloor A \rfloor + 1 = \Delta A + 1$. Since $ \lfloor A \rfloor \cup e $ is a subgraph of $A$, by Claim \ref{deltaxleqdeltaa}, $ \Delta (\lfloor A \rfloor \cup e) \leq \Delta A$, a contradiction. 
\\[1.7ex]
${\bf (p1)}$ 
We prove $(c1)$.
First we claim that  $A \setminus e$ is not connected. Suppose, to the contrary, $A \setminus e$ is connected.
Then $\Delta (A \setminus e) = \Delta A  - 1$. By Claim \ref{deltaxleqdeltaa},
$X \subseteq E(A \setminus e) \Rightarrow \Delta A\langle X \rangle
\le \Delta A  - 1$. Let $D = E \lfloor A \rfloor$. Then by Claim \ref{coreconnected} and Definition \ref{ConGraphCore},  $\Delta A\langle D \rangle = \Delta A$.  Therefore
$E \lfloor A \rfloor = D \not \subseteq  E(A \setminus e)$, and so $e \in E \lfloor A \rfloor $, 
a contradiction.
\vskip 0.5ex
Now since $A$ is connected and $A \setminus e$ is not connected, we have: $A \setminus e$ has two components $A_1$ and $A_2$, and clearly, one of them, say $A_2 = A \setminus A_1$, has a cycle. We claim that $A_1$ is a tree.
Indeed, if not, then again $\Delta (A \setminus e) = \Delta A  - 1$, and as above, $e \in E \lfloor A \rfloor $,  a contradiction.
\vskip 0.5ex
Since $A_1$ is a tree and $e$ has one end-vertex in $A_1$ and the other end-vertex in $A_2$, we have:  $\Delta A = \Delta (A \setminus A_1) =  \Delta ( A_2)$. Therefore, by Claim \ref {coreconnected} and Definition \ref{ConGraphCore},  $A_2$ contains $E \lfloor A \rfloor$. 
\\[1.7ex]
${\bf (p2)}$
We prove $(c2)$. Since $A$ has one cycle $C$, then 
$E \lfloor A \rfloor = E(C)$, and so $A \setminus e$ is a tree.
 \\[1.7ex]
${\bf (p3)}$ 
We prove $(c3)$.
If $ A \setminus e$ is connected, then obviouly $ A \setminus e$ contains a cycle.
So we assume that   $ A \setminus e$ is not connected. Then $ A \setminus e$ has two components, say $A_1$ and $A_2$ and $e$ has one end-vertex 
in $A_1$ and the other end vertex  in $A_2$,  and so $e \not \in E( A \setminus A_1)$. Clearly, at least one of $A_1$ and $A_2$, say $A_2  = A \setminus A_1$, contains a cycle. 
\vskip 0.5ex
Now, suppose to the contrary that $A_1$ has no cycle, i.e. $A_1$ is 
a tree. 
\\
Then 
$\Delta A = \Delta (A \setminus A_1) =  \Delta ( A_2)$. 
By Claim \ref{coreconnected} and Definition \ref{ConGraphCore}, the core of $A$ is an edge subset of 
$ A_2$. Therefore $e \not \in E \lfloor A \rfloor$, a contradiction.
\\[1.7ex]
${\bf (p4)}$ 
Finally, we prove $(c4)$. 
\vskip 0.5ex
First suppose that $e \notin E \lfloor A \rfloor$. Then by $(c1)$, $A \setminus e$ has a component $D$ containing $E \lfloor A \rfloor$. Since $\lfloor A \rfloor$ has at least one cycle, clearly $D$ also has a cycle. Since vertex $v$ is incident to $E \lfloor A \rfloor$, then $v$ is in $D$.
\vskip 0.5ex
Now suppose that $e \in E \lfloor A \rfloor$. Since $A \setminus e$ has two components, by $(c2)$, $A$ has at least two cycles. Then, by $(c3)$, every component of $A \setminus e$ contains a cycle. 
\ep

\vskip 1.7ex 
From Claim \ref{A-e} $(c1)$ and $(c3)$ we have: 

\begin{claim}
\label{CoreEdge,ConG}

Let $A$ be a connected graph with at least two cycles and $e \in E(A)$. Then the following are equivalent:
\\[1ex]
$(a1)$ $e \in E \lfloor A \rfloor = E [ A ]$ and
\\[1ex]
$(a2)$ every component of $A \setminus e$ has a cycle $($or, equivalently, 
$A \setminus e$ has no tree component$)$.

\end{claim}  

\vskip 0.7ex 

From Claim \ref{CoreEdge,ConG} we have the following characterization of the edges in a graph belonging to the core.

\begin{Theorem}
\label{CoreEdge,NonConG}

 Let $G$ be a graph, $t(G)$ the number of tree components of $G$, and 
$Y(G) = \{a \in E(G): t(G \setminus a) > t(G) \}$. Suppose that $G$ has a component with at least two cycles.  

 Then the following are equivalent:
\\[1ex]
$(a1)$ $e \in E [ G ] $  and
\\[1ex]
$(a2)$ $e \in E(G) \setminus Y(G)$. 

\end{Theorem}

\bp (uses Definition \ref{NonConGraphCore} and Claims \ref{A-e} and \ref{CoreEdge,ConG})
\vskip 0.5ex
If $G$ is connected, then our claim follows from Claim \ref{CoreEdge,ConG}. So we assume that $G$ is not connected. 
\vskip 0.5ex 
First, we prove $(a1) \Rightarrow (a2)$. Suppose that $(a1)$ holds, i.e. 
$e \in E [ G ] $. Then by Definition \ref{NonConGraphCore}, 
$e \in E \lfloor A \rfloor$, where $A$ is a component of $G$ with at least two cycles and $\lfloor A \rfloor = [A]$ . Then by Claim \ref{CoreEdge,ConG}, $A \setminus e$ has no tree component, and so $t(G \setminus e) = t(G)$. Thus, $e \notin Y(G)$. Therefore 
$(a2)$ holds.  
\vskip 0.5ex
Now we prove $(a2) \Rightarrow (a1)$. Suppose that $(a2)$ holds, i.e. 
$e \in E(G) \setminus Y(G)$. Then $t(G \setminus e) = t(G)$. Let $A$ be the component of $G$ containing $e$. Since $t(G \setminus e) = t(G)$, component $A$ has at least two cycles. We claim that $e \in E \lfloor A \rfloor$. Suppose, not. Then by Claim \ref{A-e} $(c1)$, $A \setminus e$ has a tree component, say $T$. Then $T$ is a component of 
$G \setminus e$. Clearly, $T$ is not a component of $G$ because $T$ is a subgraph of the component $A$ of $G$. Therefore $t(G \setminus e) > t(G)$, a contradiction. Thus, $e \in E \lfloor A \rfloor = E [A]$, and therefore 
$e \in E [ G ] $. Therefore $(a1)$ holds.        
\ep

\section{On the $k$-circular matroid of a graph}
\label{k-matroid}

\indent

In this Section we will introduce and study the properties of the so called {\em $k$-circular matroid $M_k(G)$ of a graph $G$}, where $k$ is a non-negative integer. 
We will see, in particular, that  (as before) $M_0(G)$ is the cycle matroid of graph $G$ and $M_1(G)$ is the  bicircular matroid of graph $G$.

\vskip 1ex 

The results of this Section will provide, in particular, a proper basis for our study of the problem 
$ (WP) _k$ on describing the classes of graphs with the same $k$-circular matroids.

\subsection{Circuits of the $k$-circular matroid of a graph}
\label{Circuits}

\begin{definition} {\sc Non-decreasing and submodular functions}

\vskip 0.6ex
{\em
Let $E$ be a finite set. Consider function $f: 2^{E} \rightarrow \Z$. 
\\[0.5ex]
Function $f$ is called {\em non-decreasing} if $X \subseteq Y \subseteq E \Rightarrow f(X) \le f(Y)$. 
\\[0.5ex]
Function $f$ is called {\em submodular} if $ X \subseteq Y \subseteq E \Rightarrow f(X \cup Y) + f(X \cap Y) \le f(X) + f(Y)$. 
}
\end{definition}

\begin{Theorem} {\em (J. Edmonds and  G.-C. Rota \cite{Edm},  see also \cite{Ox})}
\label{submodular}

\vskip 0.6ex
Let $E$ be a finite set and $f: 2^{E} \rightarrow \Z$ be a non-decreasing and submodular function. 
Then
$\mathcal{C}(f) = {\cal M}in \{ C\subseteq E: C \ne \emptyset ~ and ~ 
|C| > f(C)\}$ is the collection of circuits of a matroid on E. 

\end{Theorem}

Notice that if $f(X) < 0$ for every $X \subseteq E$, then 
$\mathcal{C}(f) = \emptyset $, and therefore $E$ is the only base of the matroid on $E$ induced by function $f$.

\vskip 1.5ex

The following is a simple strengthening of Theorem \ref{submodular}.

\begin{claim}
\label{eqsubmodular} 

Let E be a finite set and $f : 2^{E} \rightarrow \Z$ be a non-decreasing and submodular function. Then
 ${\cal C}(f) = {\cal M}in~ \{ C \subseteq E: C \ne \emptyset ~ and ~ 
 |C| = f(C) + 1\}$. 

\end{claim} 

\bp   
Let 
${\cal U} = \{ C \subseteq E: C \ne \emptyset  ~ and ~ |C| > f(C)\}$ and 
${\cal V} = \{ C \subseteq E: C \ne \emptyset ~ and ~ |C| = f(C) + 1\}$. 
It is clear that ${\cal V} \subseteq {\cal U}$. We want to show that 
${\cal M}in ~{\cal U} = {\cal M}in~ {\cal V}$. 
\\[1,5ex]
${\bf (p1)}$ First we prove that 
${\cal M}in ~{\cal U} \subseteq {\cal M}in ~{\cal V}$, i.e. that every minimal element of ${\cal U}$ is also a minimal  element of ${\cal V}$.
Let $X \in {\cal M}in ~{\cal U}$. 

\vskip 1ex 

First we show that  $X \in {\cal V}$.
Since $X \in {\cal M}in ~{\cal U}$, clearly  $X \in {\cal U}$, and therefore $X \ne \emptyset $. Let $e \in X$. 
By minimality of $X$  in ${\cal U}$, 
$|X| - 1 = |X \setminus e| \leq f(X \setminus e)$.
Since $f$ is a non-decreasing function, 
$ f(X \setminus e) \leq f(X)$. 
Since $X \in {\cal U}$, clearly $f(X) < |X|$. Now since in addition, $f$ is an integer-valued function,   we have: $f(X) \le |X|  - 1$.
Thus,
$|X| - 1 = |X \setminus e| \leq f(X \setminus e) \leq f(X) \le |X| - 1$. 
It follows that $|X| = f(X) + 1 $, and so $X \in {\cal V}$.

\vskip 1ex 

Now we will show that $X$ is a minimal  element of  ${\cal V}$, i.e.
that  $X \in {\cal M}in~ {\cal V}$. Suppose, not. Then there exists
$Y \subset X$ such that  $Y \in {\cal V}$. 
Since ${\cal V} \subseteq {\cal U}$, we have:  $Y \in {\cal V}$ and $Y \subset X$. Therefore $X \not \in {\cal M}in ~{\cal U}$, a contradiction.
\\[1.5ex]
${\bf (p2)}$ 
Now we will prove that 
${\cal M}in ~{\cal U} \supseteq {\cal M}in ~{\cal V}$. 
Let $Y \in {\cal M}in ~{\cal V}$. Since ${\cal V} \subseteq {\cal U}$, clearly 
$Y \in {\cal U}$. Since ${\cal U}$ is a finite family,  $Y$ must contain a subset, say $Z$,  which is a minimal  element in ${\cal U}$. However, we have proved in ${\bf (p1)}$ that ${\cal M}in ~{\cal U} \subseteq {\cal M}in ~{\cal V}$,  and therefore  $Z \in {\cal M}in ~{\cal V}$. 
It follows that 
$Y = Z \in {\cal M}in ~{\cal U}$.

\vskip 1ex 

Thus, ${\cal M}in ~{\cal U} = {\cal M}in ~{\cal V}$.        
\ep    

\vskip 2ex

Consider a function $f_{k}: 2^{E}\rightarrow \Z$  such that
$f_{k}(X) = |V(G\langle X \rangle )| - 1 + k$ for every  $X \in  2^{E}$. 
It is easy to see that the  following is true.

\begin{claim} {\em (\cite{Edm}, see also \cite{Ox} )}
\label{gsubmodular}
\vskip 0.5ex
Let   
$G =(V, E,\phi)$
be a graph and $k \ge 0$. 
Then $f_{k}$ is a non-decreasing and submodular function.
\end{claim}

\vskip 1ex

From Claims \ref{eqsubmodular} and   \ref{gsubmodular} we have:
\begin{claim}
\label{Ck(G)=} 
Let $G = (V, E, \phi )$ be a graph and 
$k \ge 0$.
Then 
\vskip 0.7ex
${\cal C}_{k}(G) = {\cal M}in ~\{ C \subseteq E: C \ne \emptyset ~and 
~|C| = | V(G\langle C \rangle )| + k \}$  
\\[0.7ex]
is the collection of circuits of a matroid on $E$.
\end{claim} 

\begin {definition} {\sc The $k$-circular matroid of a graph}
\vskip 0.7ex

{\em
We call the matroid described in Claim \ref {Ck(G)=} by the set ${\cal C}_{k}(G)$ of its circuits 
{\em the $k$-circular matroid of} $G$ and  denote it $M_k(G)$.
Let ${\cal D}_k(G)$, ${\cal I}_k(G)$, ${\cal B}_k(G)$, and ${\cal C}^*_k(G)$ denote the families of dependent sets, independent sets, bases, and cocircuits of $M_k(G)$, respectively. 
}
\end{definition}

It is easy to see  that (as above)  $M_0(G)$ is the cycle matroid  of a graph $G$.
As we will see below,  matroid $M_1(G)$ is  the  bicircular matroid $B(G)$.  

\vskip 2ex

From Claim \ref{Ck(G)=}  we have: 
\begin{claim}
\label{deltageqk}
Let $G = (V, E, \phi )$ be a graph and $k \ge 0$. Then 
\vskip 1ex 

$\mathcal{C}_{k}(G)=
{\cal M}in~\{C\subseteq E: C \ne \emptyset~ and ~\Delta G\langle C\rangle  = k \} $. 
\end{claim}

As in Section 1, ${\cal G}_{\bowtie}$ is the set of graphs with no isolated vertices, no leaves, and  no cycle components, i.e. the set of cacti-graphs.

\begin{Theorem} 
{\sc Structure of graph induced by a circuit of $ M_{k}(G)$ in $G$}
\label{cirstructure}
\vskip 0.7ex
Let $G = (V, E,\phi)$ be a graph. Let $C \subseteq E$ and
 $k \ge 1$. 
 Then the following are equivalent:
\\[1ex]
$(c1)$ $C \in \mathcal{C}_k(G)$ and
\\[1ex]
$(c2)$ $\Delta G \langle C \rangle = k$ and 
$G \langle C \rangle \in {\cal G}_{\bowtie}$, i.e. 
$G\langle C \rangle$ has no isolated vertices, no leaves, and every component of $G\langle C \rangle$ has at least two cycles. 

\end{Theorem}

\bp (uses Claim \ref{deltageqk})
\vskip 0.8ex
Obviously, a connected graph $A$ has at least two cycles if and only if 
$\Delta A \ge 1$. 
\\[1ex] 
${\bf (p1)}$ First we prove $(c1) \Rightarrow (c2)$. 
By Claim \ref{deltageqk}, $\Delta G \langle C \rangle = k$. 
Also, if $e$ is a pendant edge of $G \langle C \rangle$, then clearly, 
$\Delta G \langle C \rangle - \Delta G \langle C\setminus e \rangle \in \{-1,0\}$. Therefore 
$\Delta G \langle C\setminus e \rangle \geq \Delta G \langle C \rangle$, 
contradicting minimality of $C$. Hence  $G \langle C \rangle$ contains no leaves. 

\vskip 1.2ex 

Let ${\cal A}$ be the set of components of $G \langle C \rangle$. Clearly, 
$E(A) \neq \emptyset$ for every  
$A \in {\cal A}$,  
\\[0.8ex]
$C = \bigcup_{A\in{\cal A}} E(A)$, and $V(G \langle C \rangle) = 
\bigcup_{A\in{\cal A}} V(A)$. Then
\\[1.2ex]
$~~~~~~~~~~~~\sum_{A\in{\cal A}} \Delta A  = 
\sum_{A\in{\cal A}} (|E(A)| - |V(A)|) =
\sum_{A\in{\cal A}} |E(A)| - \sum_{A\in{\cal A}} |V(A)|
\\[0.8ex]
~~~~~~~~~~~~~~~~~~~~~~~~~~= |C| - |V(G \langle C \rangle)| = \Delta G \langle C \rangle = k \geq 1$. 
\\[1.2ex]
 Therefore, 
 $\sum_{A\in{\cal A}}   \Delta A
  = \Delta G \langle C \rangle \geq 1$.  

\vskip 1ex 

If $|{\cal A}| = 1$, then
 $\Delta A \geq 1$, where $A$ is the unique component of $G[C]$. Now 
 suppose that $|{\cal A}| \geq 2$ and 
let $A' \in {\cal A}$. Then the set of components of $G  \langle C \setminus E(A') \rangle$
\\[0.8ex]
 is 
${\cal A}' = {\cal A} \setminus \{A'\}$. 
Therefore
\vskip 0.8ex
$\Delta G\langle C \setminus E(A') \rangle  = 
\sum_{A \in {\cal A}'} \Delta A =
\sum_{A \in {\cal A}}\Delta A - \Delta A' =
\Delta G  \langle C \rangle - \Delta A'$. 

Thus, if  
$\Delta A' < 1$,
then $\Delta G  \langle C \setminus E(A') \rangle \geq \Delta G \langle C \rangle = k$, contradicting
 minimality of $C$. Hence every component $A$ of 
$G\langle C \rangle$ has at least two cycles. 
\\[2ex]
${\bf (p2)}$ Finally,  we prove  $(c1) \Leftarrow (c2)$. Consider $C \subseteq E$ satisfying $(c2)$. Then   $\Delta G \langle C  \rangle = k$, where $k \ge 1$. Therefore $C \ne \emptyset $. Hence 
$C$ is an element of set ${\cal V} = \{Z\subseteq E: Z \ne \emptyset~ and ~
\Delta G \langle Z \rangle = k\} $. It remains to show that $C$ is a minimal element in  ${\cal V}$.
%. 

\vskip 1ex

We first prove the following
\vskip 1ex
{\sc Claim.}
{\em 
If $A$ is a component of $G \langle C  \rangle $
and $X \subset E(A)$, 
then $\Delta G\langle X\rangle < 
\Delta A$.
}
\\[1ex]
{\em Proof.}
If $X = \emptyset$, then $\Delta G\langle X \rangle = 0 < 1 \leq 
\Delta A $.
So we assume that $X \ne \emptyset $.
\vskip 0.3ex
Let $S = V(A)\setminus V(G\langle X \rangle )$.
 If $S =\emptyset$, then 
$V(G\langle X \rangle ) = V(A)$ and since $|X| < |E(A)|$, 
\\[0.3ex]
we have:
 $\Delta G\langle X \rangle =| X| - |V(G\langle X \rangle )| < |E(A)| - |V(A)| = 
 \Delta A $.
 So we assume that $S \ne \emptyset$. 
 \\[0.3ex]
 Since $X \ne \emptyset$ and $A$ is connected and has no leaves, 
 there are at least $|S |+ 1$ edges of $E(A) \setminus X$ 
 \\[0.3ex]
 that are incident to vertices in 
 $S$,  
 and so 
 $|E(A)\setminus X| > |S|$. 
 Then 
 \vskip 0.7ex
~~~~~~~$\Delta G \langle E(A)\rangle = E(A ) - V(A) = 
 |X| + |E(A) \setminus X| - (|V(G\langle X \rangle)| + |S|) 
  \\[0.7ex]
 ~~~~~~~~~~~~~~~~~~~~~~~~~~~= 
 \Delta G\langle X \rangle + |E(A)\setminus X| - |S| \geq \Delta G \langle X \rangle + 1 > \Delta G\langle X \rangle $. 
\epcl
\vskip 1ex
Now we prove that $C$ is a minimal element in  ${\cal V}$, namely, that 
\\
$Z \subset C ~ and~ Z \ne \emptyset \Rightarrow
Z \not \in {\cal V}$, i.e.  
$\Delta_{G}(Z) < k$. 

\vskip 1ex 

Let ${\cal A}$ be the set of components of $G \langle  C \rangle$. 
Put $Z_A = Z \cap E(A)$ for $A \in {\cal A}$.
Then 
\\[1.8ex]
$\Delta G\langle Z \rangle  = |Z| - |V(G\langle Z \rangle )| = 
\sum_{A \in{\cal A}}  |Z \cap E(A)| - 
\sum_{A \in {\cal A}}| V(G\langle Z \rangle )\cap V(A)| =
\\[1ex]
\sum_{A \in {\cal A}}(|Z\cap ~E(A)| - |V(G\langle Z \rangle ) \cap ~V(A)|) =
\sum_{A \in {\cal A}} \Delta G\langle  Z_A\rangle $. 
\\[1.8ex]
Since $Z \subset Y$, we have: $Z_A \subset A$ for some $A \in {\cal A}$. Now by the above {\sc Claim}, 
\\[1.6ex]
$\Delta G\langle Z \rangle = \sum_{A\in{\cal A}}\Delta G\langle  Z_A\rangle  < ~\sum_{A \in {\cal A}}
\Delta A = 
\Delta G \langle  C \rangle = k$.  
\ep 

\vskip 2ex

We remind that a graph $R$ is a bicycle if and only if $R$ is a connected graph with no leaves and $\Delta R = 1$ (see Figure \ref{Bicycles2Fig}).

\vskip 2ex
From the definition of a bicycle, Claim \ref{BicyMin}, and Theorem \ref{cirstructure} we have: 

\begin{claim}
\label{BicyC_1}
 
$C \in \mathcal{C}_1(G)$ if and only if $G \langle C \rangle$ is a bicycle {\em (see Figure \ref{Bicycles2Fig})}. 
 
\end{claim}

\begin{remark}
\label{ConAndNonconCircuits}
{\em
There is an essential difference between bicircular matroid $M_1(G)$ and 
\\
$k$-circular matroid $M_k(G)$ for $k \ge 2$. Namely, every circuit of $M_1(G)$ induces a connected subgraph in $G$. Therefore if $M_1(G)$ is a connected matroid, then $G$ is a connected graph. On the other hand, if $k\ge 2$, then $M_k(G)$ may have a circuit that induces a non-connected subgraph in $G$. Therefore for $k\ge 2$ matroid $M_k(G)$ may be connected although graph $G$ is not connected.   
}
\end{remark}

\vskip 3ex
A circuit $C$ of $M_k(G)$ is called {\em graph-connected} if the graph induced by $C$ in $G$ is connected. 

\vskip 2ex
Matroid $M_2(G)$ is called a {\em tricircular matroid}. The connected circuits of tricircular matroids are subdivisions of graphs in Figure \ref{TricyclesFig}. The non-connected circuits of tricircular matroids are pairs of disjoint subdivisions of graphs from $\{S, D, B\}$ in Figure \ref{Bicycles2Fig}.

\vskip 4ex

\begin{figure}[h]
\begin{center}
\scalebox{0.28}[.28]{\includegraphics{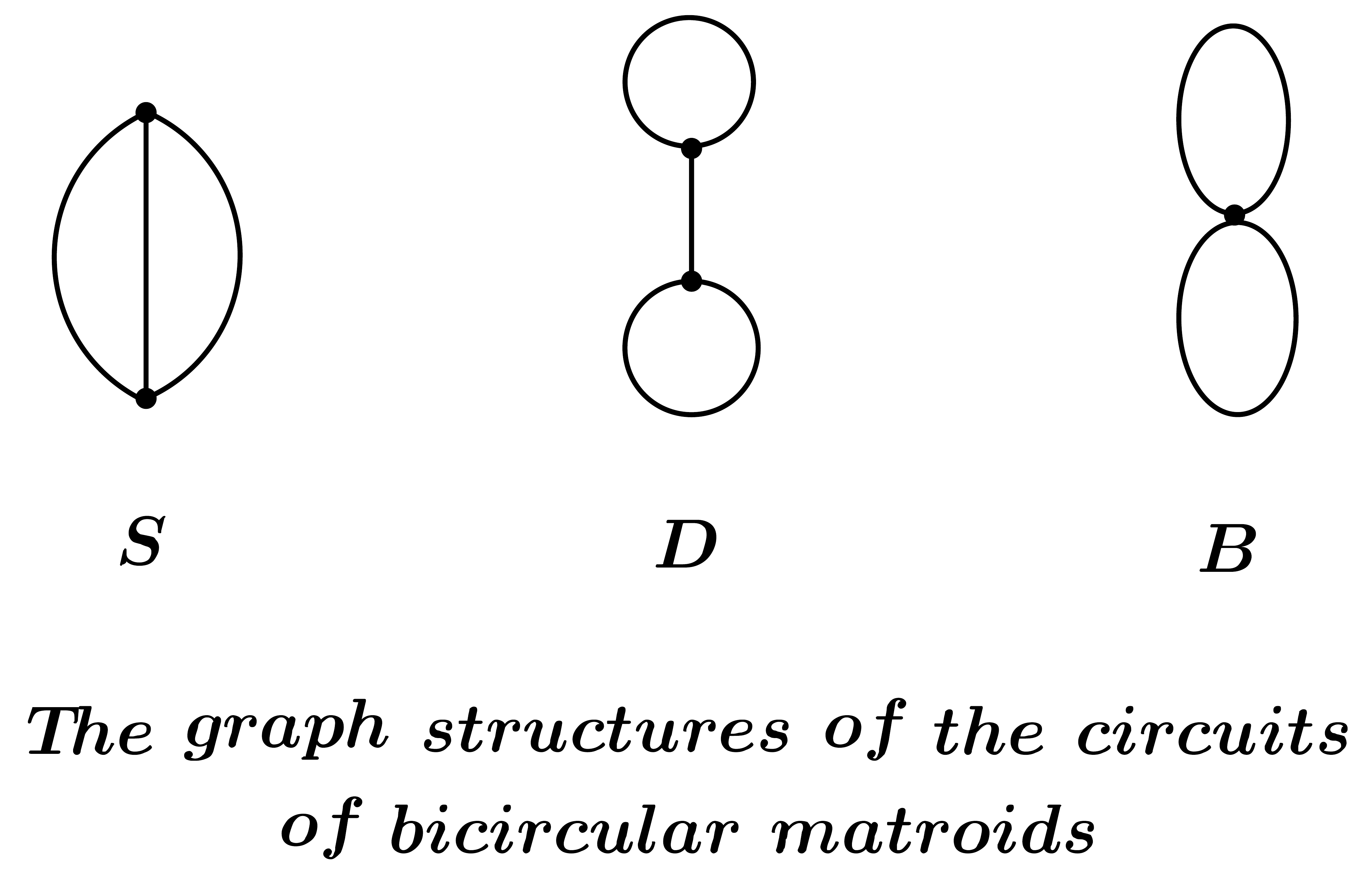}}
\end{center}
\caption
{Bicycles (up to subdivisions).}
\label{Bicycles2Fig}
\end{figure}

\begin{figure}[h]
\begin{center}
\scalebox{0.29}[.29]{\includegraphics{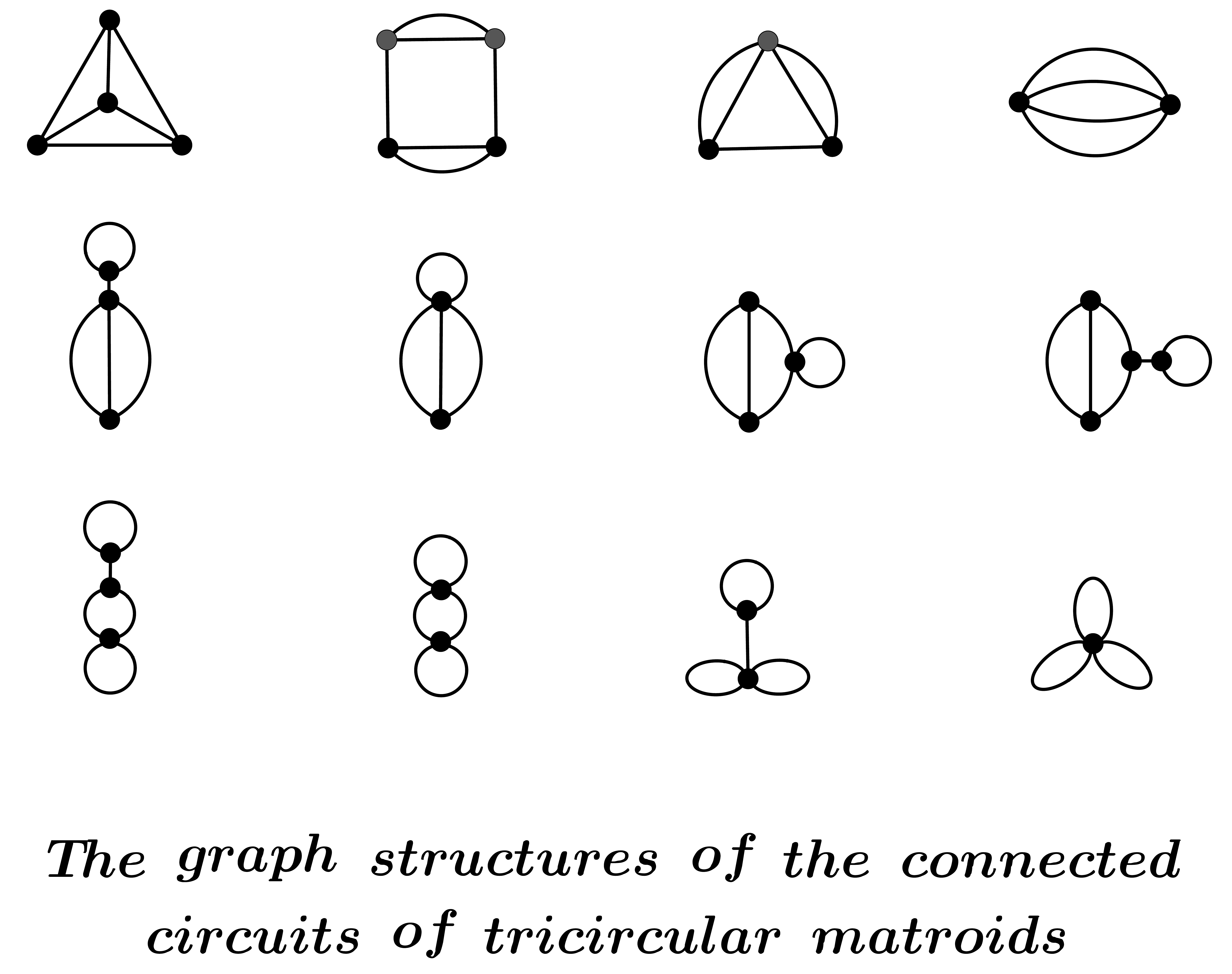}}
\end{center}
\caption
{Connected tricycles (up to subdivisions).}
\label{TricyclesFig}
\end{figure}

\clearpage

\begin{Theorem} 
{\sc Recursive description of a circuit of $ M_{k}(G)$}
\label{CirRecDesc}
\vskip 0,7ex
Let $G = (V, E,\phi)$ be a graph. Let $C \subseteq E$ and
 $k \ge 1$. If $k = 1$, then $C \in \mathcal{C}_1(G)$ if and only if $G \langle C \rangle$ is a bicycle. If $k \ge 2$, then the following are equivalent:
\\[1.2ex]
$(c1)$ $C \in \mathcal{C}_k(G)$ and
\\[1,2ex]
$(c2)$ there exists $C' \in {\cal C}_{k-1}(G)$ such that 
$G \langle C \rangle = G \langle C' \rangle \cup P$, where $P$ is either a path or a cycle or a lollipop {\em(see Figure \ref{Ears})} or a bicycle and  
\\[1.2ex]
${\bf (l,p)}$ $V(G \langle C' \rangle) \cap V(P) = End(P)$ if $P$ is a lollipop or a path,
\\[1.2ex] 
${\bf (c)}$ $|V(G \langle C' \rangle) \cap V(P)| = 1$ if $P$ is a cycle, and 
\\[1,2ex]  
${\bf (b)}$ $V(G \langle C' \rangle) \cap V(P) = \emptyset$ if $P$ is a bicycle. 

\end{Theorem}

\bp (uses Theorems \ref{EarDcmpGinfty} and \ref{cirstructure} and Claim \ref{a,bEarDcmpGinfty})
 \vskip 0.7ex
If $k = 1$, then by Claim \ref{BicyC_1}, $C \in \mathcal{C}_1(G)$ if and only if $G \langle C \rangle$ is a bicycle. So we assume that $k \ge 2$. 
 \vskip 0.7ex
First, we prove $(c1) \Rightarrow (c2)$. Let $C \in \mathcal{C}_k(G)$. Then by Theorem \ref{cirstructure}, $\Delta G \langle C \rangle = k$ and 
$G \langle C \rangle  \in {\cal G}_{\bowtie}$. By Claim \ref{a,bEarDcmpGinfty}, graph 
$G \langle C \rangle$ contains a bicycle $D$. Therefore by Theorem \ref{EarDcmpGinfty}, there exists a subgraph $G_{k-1} = G \langle C' \rangle$ such that $C' \in {\cal C}_{k-1}(G)$ and $(c2)$ holds.  

\vskip 0.7ex
Now we prove $(c2) \Rightarrow (c1)$. Since $(c2)$ holds, $G \langle C \rangle  \in {\cal G}_{\bowtie}$ and 
$\Delta G \langle C \rangle = \Delta G \langle C' \rangle + 1 = k$. Therefore by Theorem \ref{cirstructure}, $C \in \mathcal{C}_k(G)$, and so $(c1)$ holds.  
\ep

\begin{definition}
\label{Mk-Calsses}
{\em
Given two graphs $G$ and $G'$, we say that $G$ and $G'$ are 
{\em $M_k$-equal } and write $G \sim^k G'$ if $M_k(G) = M_k(G')$. Clearly, relation $\sim^k$ on the set of finite graphs is an equivalence relation. The equivalence classes of this relation are called {\em $M_k$-equivalence classes}. 
}
\end{definition}

Using Theorems \ref{EarDcmpGinfty} and \ref{cirstructure} it is not difficult to prove the following claim about the relation between $M_k$- and $M_s$-equivalence classes of graphs for $k > s \ge 1$. 

\begin{claim}
\label{Mk,sEquivalence}  
Let $k > s \ge 1$. If ${\cal F}$ be an $M_k$-equivalence class of graphs, then ${\cal F}$ is the union of some $M_s$-equivalence classes of graphs. In other words, if $M_s(G) = M_s(G')$, then $M_k(G) = M_k(G')$. 
\end{claim}

\begin{claim}
\label{no-loops}

Let $G$ be a graph and $k \geq 1$. Then $M_k(G)$ has no loops.

\end{claim}

\bp (uses Theorem \ref{cirstructure})
 \vskip 0.3ex
By Theorem \ref{cirstructure}, if $C$ is a circuit of $M_k(G)$, then 
$G \langle C \rangle$ has at least two cycles. Therefore $|E(C)| \ge 2$. 
\ep

\begin{claim}
\label{geqkdep}

Let $G = (V, E,\phi )$ and $k \geq 1$. If $X \subseteq E$ and $\Delta G\langle  X \rangle \geq k$, then $X \in \mathcal{D}_k(G)$.

\end{claim}

\bp
By Claim \ref{eachi}, the set 
${\cal P} = \{Z\subseteq X: \Delta G\langle  Z \rangle = k\}$ is non-empty. 
A minimal 
\\[0.9ex]
element of ${\cal P}$ exists and belongs to ${\cal C}_{k}(G)$, implying that $X\in {\cal D}_{k}(G)$.  
\ep

\subsection{Non-trivial $k$-circular matroids}
\label{Non-trivialMk}

\indent

We remind that a matroid $M = (E, {\cal I})$ is non-trivial if $E$ is not a base and not a cobase of $M$.

In this Subsection we give a criterion for a $k$-circular matroid of a graph to be non-trivial and, on the other hand, describe all graph representations of a trivial $k$-circular matroid of a graph. 

\vskip 1ex

\begin{Theorem} {\sc  A criterion for matroid $M_{k}(G)$ to be non-trivial}
\label{maxk}

Let $G$ be a graph, $F(G)$ the union of all tree components of $G$,
 and $k \geq 1$. Then the following are equivalent:    
\\[0.7ex]
$(a1)$ ${\cal C}_k(G) \ne \emptyset$,       
\\[0.5ex]
$(a2)$ $k \leq \Delta G + cmp(F(G))$, and 
\\[0.5ex]
$(a3)$ $M_k(G)$ is a non-trivial matroid.     

\end{Theorem} 

\bp  (uses Theorem \ref{cirstructure} and Claims \ref{deltaxleqdeltah} and \ref{geqkdep})
\vskip 0.5ex
Let $G' = G \setminus F(G)$. Clearly,  
$ \Delta G = \Delta G' + \Delta F(G) = \Delta G' - cmp(F(G))$. Therefore 
\\[0.3ex] 
$\Delta G' = \Delta G + cmp(F(G))$. 

\vskip 1ex

First we prove $(a1) \Rightarrow (a2)$. Suppose that $(a1)$ is true, i.e. ${\cal C}_{k}(G) \ne \emptyset$. Let $C \in {\cal C}_{k}(G)$. Then by Theorem \ref{cirstructure}, $\Delta G \langle C \rangle \ge k$ and every component of $G \langle C \rangle$ has at least two cycles. Since every vertex of $F(G)$ belongs to a tree components of $G$, graphs $G \langle C \rangle$ and 
$F(G)$ have no vertex in common. Since $G' = G \setminus F(G)$, graph 
$G \langle C \rangle$ is a subgraph of  $G'$. Since $G'$ has no tree components, by Claim \ref{deltaxleqdeltah},  
$k \le \Delta G \langle C \rangle \le \Delta G' = \Delta G + cmp(F(G))$.  
 
 \vskip 1ex    

Next we prove $(a2) \Rightarrow (a1)$. Suppose that $(a1)$ is true, i.e. 
$k \leq \Delta G + cmp(F(G))$. Since $\Delta G' = \Delta G + cmp(F(G))$, we have $\Delta G' \ge k$. Then by Claim \ref{geqkdep}, 
$E(G') \in \mathcal{D}_k(G)$ and therefore ${\cal C}_{k}(G) \ne \emptyset$. 

 \vskip 1.5ex 
 
Now we prove $(a1) \Leftrightarrow (a3)$.
By definition of a non-trivial matroid, $(a3) \Rightarrow (a1)$. By Claim 
\ref{no-loops}, $M_k(G)$ has no loops, and therefore $E(G)$ is not a cobase of $M_k(G)$. By $(a1)$, $E(G)$ is not a base of $M_k(G)$. Thus, $M_k(G)$ is a non-trivial matroid, and therefore $(a1) \Rightarrow (a3)$.
\ep
  \vskip 1ex
 
From Theorem \ref{cirstructure} and Claim \ref{maxk} we have:

\begin{Theorem}
\label{GRepTriv}

 Let $G$ be a graph. Then $M_k(G)$ is a trivial matroid if and only if $G$ has no subgraph $F$ such that $F \in {\cal G}_{\bowtie}$ and $\Delta F = k$. 

\end{Theorem}

\begin{remark}
\label{ReducToNonTrivialM}
{\em
The above Theorem \ref{GRepTriv} provides for a given $k \ge 1$ a complete description of all graphs $G$ that have the same matroid $M_k(G)$ in case when $M_k(G)$ is a trivial matroid. In other words, the above Theorem \ref{GRepTriv} describes all graph representations of the $k$-circular trivial matroid. 
}
\end{remark}

For this reason from now on we will consider the situations when $M_k(G)$ is a non-trivial matroid.

\subsection{Coloops of non-trivial  $k$-circular matroids} 
\label{MkColoops}

\indent

Now we are going to give a graph characterization of the set of coloops for a non-trivial matroid $M_k(G)$. 

 \vskip 1ex
 
We know that an element $e$ in a matroid $M$ is a coloop of $M$ if and only if $e$ belongs to no circuit of $M$. 

\vskip 1ex 

Obviously, from Claim \ref{maxk} we have:  

\begin{claim}
\label{}
 
$M_k(G)$ is a non-trivial matroid if and only there exists an edge in $G$ which is not a coloop of $M_k(G)$, or, equivalently, $M_k(G)$ is a trivial matroid if and only if every edge of $G$ is a coloop of $M_k(G)$.  
 
\end{claim}

\begin{claim} 
\label{tree-components}
 
Let $G$ be a graph and $R = R(G)$ the union of
 components of $G$ having at most one cycle.
 Then every edge $e$ of $R$ is a coloop of matroid $M_k(G)$
 for every positive integer $k$, and so $E(R)$ is a subset of every base of $M_k(G)$.  
 
\end{claim}

\bp  (uses Claim \ref{cirstructure})
\vskip 0.5ex
Obviously, an element $e$ of a matroid $M$ is a coloop if an only if belongs to no circuit of $M$. Let $e \in E(R)$, i.e. $e$ belongs to a component of $G$, say $T$, having at most one cycle.

 \vskip 1ex
 
It is sufficient to prove that $e $ belongs to no circuit of $M_k(G)$.
Suppose, to the contrary, that there exists a circuit $C$ in $M_k(G)$ containing $e$. Let $A$ be the component of $G\langle C \rangle $ containing $e$.
Then by Claim \ref{cirstructure}, $A $ has at least two cycles $D_1$ and $D_2$. Now since $e \in E(A) \cap E(T)$, clearly $A$ is a subgraph the  component $T$, and so cycles $D_1$ and $D_2$ are subgraphs of $T$. 
\ep

 \vskip 2ex
 
Notice that if $M_k(G)$ is a non-trivial matroid for some $k \ge 1$, then by Claim \ref{maxk}, there exists $ C \in {\cal C}_{k}(G) $. By Theorem \ref{cirstructure}, every component of $ G \langle C \rangle $ has at least two cycles. Therefore graph $G$ has a component with at least two cycles, and so the core $[ G ]$ of graph $G$ (see Definition \ref{NonConGraphCore}) is defined.

\vskip 1.5ex 
Let $L_k^*(G)$ denote the set of coloops of $M_k(G)$. 
 
 \begin{claim}
 \label{coloops}

Let $k \ge 1$ and $E = E(G)$. Suppose that $M_k(G)$ is a non-trivial matroid. Then $L_k^*(G) = E \setminus E [ G ]$. 

  \end{claim}
 
 \bp (uses Theorem \ref{EarDcmpGinfty}, \ref{CoreInGinfty}, and \ref{cirstructure}, Corollary \ref{a,bInFsbgrG}, and Claims \ref{minmax}, 
 \ref{A-e},  and \ref{tree-components})  
 \\[1ex]
${\bf (p1)}$ First we prove that 
$E \setminus E [ G ] \subseteq L_k^*(G)$. Let $e \in E \setminus E [ G ]$ and 
$A$ be the component of $G$ containing $e$. If $A$ has at most one cycle, then by Claim \ref{tree-components}, $e \in L_k^*(G)$. So we assume that 
$A$ has at least two cycles. Since $e \in E \setminus E [ G ]$, clearly 
$e \notin E [A] = E \lfloor A \rfloor$. Then, by Claim \ref{A-e} $(c1)$, graph 
$A \setminus e$ has a tree component. We claim that $e$ belongs to no circuit of $M_k(G)$. Suppose, to the contrary, that there exists a circuit $C$ of $M_k(G)$ containing $e$. Let $D$ be the component of 
$G \langle C \rangle$ containing $e$. By Theorem \ref{cirstructure}, $D$ has no leaf. Then by Claim \ref{minmax}, $e \in E(D) \subseteq E \lfloor A \rfloor$, a contradiction. Therefore $e$ belongs to no circuit of $M_k(G)$, and so 
$e \in L_k^*(G)$. Hence, 
$E \setminus E [ G ] \subseteq L_k^*(G)$. 
 \\[1.5ex]
${\bf (p2)}$ Now we prove that 
$L_k^*(G) \subseteq E \setminus E [ G ]$. Suppose, to the contrary, that there exists $e \in L_k^*(G)$ such that 
$e \notin E \setminus E [ G ]$, i.e. $e \in E [ G ]$. We will show that there exists $C \in {\cal C}_k(G)$ containing $e$, contradicting $e \in L_k^*(G)$. 

 \vskip 1ex
 
Since $M_k(G)$ is a non-trivial matroid, there exists $D \in {\cal C}_k(G)$. By Theorem \ref{cirstructure}, $\Delta G \langle D \rangle = k \ge 1$ and 
$G \langle D \rangle \in {\cal G}_{\bowtie}$. Therefore $G$ has a component with at least two cycles. Then by Theorem \ref{CoreInGinfty}, $[ G ] \in {\cal G}_{\bowtie}$ and $G \langle D \rangle \subseteq^s [ G ]$. Also, by Theorem \ref{EarDcmpGinfty}, $\Delta [ G ] \ge \Delta G \langle D\rangle = k$. Now by Corollary \ref{a,bInFsbgrG}, there exists a subgraph $F$ of $[ G ]$ such that 
$e \in E(F)$, $F \in  {\cal G}_{\bowtie}$, and  $\Delta (F) = k$. Put $E(F) = C$. Then $e \in C$ and by Theorem \ref{cirstructure}, $C  \in {\cal C}_k(G)$.  
\ep
 
\vskip 1.5ex 
 
Given a graph $G$, let $t(G)$ be the number of tree components of $G$, and  $Y(G) = \{e \in E(G): t(G \setminus e) > t(G) \}$.

\vskip 1ex 
From Claims \ref{CoreEdge,NonConG} and \ref{coloops} we have the following more detailed characterization of the set of coloops of a non-trivial matroid $M_k(G)$.

\begin{Theorem} 
{\sc Graph description of coloops of $ M_{k}(G)$}
\label{ColoopsCriterion}

\vskip 0.7ex
Let $G$ be a graph and $M_k(G)$ a non-trivial matroid. Then 
$L_k^*(G) = E \setminus E [ G ] = Y(G)$.

\end{Theorem}

\subsection{Connected $k$-circular matroids}
\label{MkCircConn}

\indent

We recall that a matroid $M = (E, {\cal I})$ is connected if  $|E|\ge 2$ and 
for every $a, b \in E$ there exists $C \in {\cal C}(M)$ such that $a, b \in C$.

In this Subsection we present a graph criterion for a $k$-circular matroid to be connected.

 \vskip 1.5ex
  
We start with the following claims. 

\begin{claim} 
\label{CoreCircuits}

Let $G$ be a graph for which $M_k(G)$ is a non-trivial matroid. Then 
${\cal C}_k(G) = {\cal C}_k [ G ]$.

\end{claim}

\bp  (uses Claims \ref{MatroidColoops} and \ref{coloops})

\vskip 0.6ex

Obviously, ${\cal C}_k [ G ] \subseteq {\cal C}_k(G)$. We prove 
${\cal C}_k(G) \subseteq {\cal C}_k [ G ]$. Let $C \in {\cal C}_k(G)$. Then by
\\[0.6ex]
 Claim \ref{MatroidColoops}, $C \cap L_k^*(G) = \emptyset$. By Claim \ref{coloops}, 
$C \cap (E \setminus E [ G ]) = \emptyset$. 
Therefore   $C \subseteq  E [ G ]$.   
\ep 
 
\vskip 1.5ex

From Claim \ref{CoreCircuits} we have:

\begin{claim} 
\label{CoreMatorid}

Let $G$ and $G'$ be graphs and $k \ge 1$. Suppose that $M_k(G)$ is a non-trivial matroid. Then $M_k(G) = M_k(G')$ if and only if $E(G) = E(G')$ and 
$M_k [ G ] = M_k [ G' ] $.  
 
\end{claim}

\begin{claim}
\label{CoreExistence}

If $M_k(G)$ is a non-trivial matroid, then the core $ [ G ] $ of $G$ is defined. 

\end{claim}

\bp (uses Theorems \ref{a,bEarDcmpGinfty} and \ref{cirstructure})

\vskip 0.7ex
Since $M_k(G)$ is non-trivial, there exists $C \in {\cal C}_k(G)$. By Theorem \ref{cirstructure}, $G \langle C \rangle \in {\cal G}_{\bowtie}$. By Theorem \ref{a,bEarDcmpGinfty}, $G \langle C \rangle$ contains a bicycle, and so $G$ has a component with at least two cycles. Therefore the core $[ G ]$ of $G$ is defined. 
\ep

\begin{remark}
\label{ReducToClean}
{\em
Claim \ref{CoreMatorid} and Theorem \ref{CoreInGinfty} allows us to reduce our problem $ (WP) _k$ of describing the classes of graphs with the same non-trivial $k$-circular matroid to the corresponding problem for the graphs in ${\cal G}_{\bowtie}$, i.e. for the class of cacti-graphs. 
}
\end{remark}

For this reasons from now on we will consider the situations when $G$ is a cacti-graph.  

\vskip 1.5ex

Let $L_k(G)$ and $L_k^*(G)$ denote the set of loops and coloops of $M_k(G)$, respectively. By Claim 
\ref{no-loops}, $L_k(G) = \emptyset $ and by Claim \ref{coloops}, $L_k^*(G) = E \setminus E [ G ] $ for any non-trivial $M_k(G)$. 
Therefore 
$E \setminus (L_k(G) \cup L_k^*(G)) = E \setminus L_k^*(G) = E [ G ] $. In particular, for cacti-graphs 
$L_k(G) \cup L_k^*(G)) = \emptyset$.

\vskip 2ex

In what follows we give a characterization of a  non-trivial  matroid $M_k(G)$ for a cacti-graph $G$.

\begin{Theorem} 
{\sc Graph description of connected matroid $M_k(G)$}
\label{MkCon}

\vskip 0.7ex 
\indent
Let $k \ge 2$. Then the following are equivalent:
\\[1ex]
$(a1)$ $M_k(G)$ is a non-trivial matroid and $G \in {\cal G}_{\bowtie}$ and 
\\[1ex]
$(a2)$ $M_k(G)$ is a connected matroid. 

\end{Theorem}

\bp (uses Theorems \ref{EarDcmpGinfty}, \ref{CoreInGinfty}, and \ref{cirstructure}, Corollary \ref{a,bInFsbgrG}, and Claims \ref{coloops} and \ref{CoreExistence})
\\[1.5ex]
${\bf (p1)}$ First we prove $(a1) \Rightarrow (a2)$.

Let $a, b \in E(G)$. Since $M_k(G)$ is a non-trivial matroid, there exists 
$D \in {\cal C}_k(G)$. By Theorem \ref{cirstructure}, 
$\Delta G \langle D \rangle = k$ and 
$G \langle D \rangle \in {\cal G}_{\bowtie}$. Since $G \in {\cal G}_{\bowtie}$ and 
$G \langle D \rangle \subseteq^s G$, by Theorem \ref{EarDcmpGinfty}, $\Delta G \ge \Delta G \langle D\rangle = k$. By Corollary \ref{a,bInFsbgrG}, there exists a subgraph $F$ of $G$ such that $a, b \in E(F)$, $F \in  {\cal G}_{\bowtie}$, and  $\Delta (F) = k$. Put 
$E(F) = C$. Then $a, b \in C$ and by Theorem \ref{cirstructure}, 
$C  \in {\cal C}_k(G)$. Therefore $M_k(G)$ is a connected matroid.  
\\[1.5ex]
${\bf (p2)}$ Finally, we prove $(a2) \Rightarrow (a1)$.
Since $M_k(G)$ is a connected matroid, we have: $M_k(G)$ has a $k$-circuit 
$C$, and therefore $M_k(G)$ is non-trivial. By Claim \ref{CoreExistence}, the core 
$[ G ]$ of $G$ is defined. Since $M_k(G)$ is a connected matroid, 
$L_k^*(G) = \emptyset$. Then by Claim \ref{coloops}, $E = E [ G ]$, and so 
$G = [ G ]$. Now by Theorem \ref{CoreInGinfty}, 
$G \in {\cal G}_{\bowtie}$.  
\ep

\begin{Theorem} 
{\sc Graph description of the $\sim$- equivalence classes in $M_1(G)$}
\label{Mk=1Classes}
   
\vskip 0.7ex   
Let $G$ be a graph. Then the following are equivalent:    
\\[1ex]
$(a1)$$M_1(G)$ is a non-trivial matroid and $G \in {\cal G}_{\bowtie}$ and 
\\[1ex]
$(a2)$ $E(A)$ is an $\sim$-equivalence class of $M_1(G)$ for every component $A$ of $G$. 
   
\end{Theorem}

\bp (uses Theorem \ref{a,bEarDcmpGinfty} and \ref{CoreInGinfty} and Claim \ref{coloops})
\\[1.5ex]
${\bf (p1)}$ First we prove $(a1) \Rightarrow (a2)$. Let $A$ be a component of $G$ and $a, b \in E(G)$. 

First suppose that $a$ and $b$ both belong to $A$. Since $G \in {\cal G}_{\bowtie}$, clearly also $A \in {\cal G}_{\bowtie}$. By Theorem \ref{a,bEarDcmpGinfty}, there exists a bicycle $A_0$ such that $A_0 \subseteq^s A $ and $a, b \in A_0$. Since $A_0$ is a bicycle, 
$A_0 \in {\cal C}_1(G)$, and therefore $a \sim b$. 

 \vskip 1ex
 
Now suppose that $a$ and $b$ do not belong to the same component of $G$. Then no $C \in \mathcal{C}_1(G)$ contains both $a$ and $b$ because $G \langle C \rangle$ is a bicycle, which is a connected graph.   

Therefore $E (A)$ is an $\sim$-equivalence class of $M_1(G)$. 
\\[1.5ex]
${\bf (p2)}$ Now we prove $(a2) \Rightarrow (a1)$. Let $A$ be a component of $G$.

 Since $E (A )$ is an $\sim$-equivalence class of $M_1(G)$, clearly ${\cal C}_1(G) \ne \emptyset$, and therefore matroid $M_k(G)$ is non-trivial and the core $ [G] $ of $G$ is defined.    

Since $E(A)$ is an $\sim$-equivalence class of $M_1(G)$ for every component $A$ of $G$, we also have: every component $A$ has at least two elements that belong to a common circuit of $M_1(G)$. Therefore $L_k^*(G) = \emptyset$. Then by Claim \ref{coloops}, $E = E [G]$, and so $G = [G]$. Now by Theorem \ref{CoreInGinfty}, 
$G \in {\cal G}_{\bowtie}$.  
\ep

\begin{corollary}
\label{CmpM_1(G)}

Let $G$ be a graph. Suppose that $M_1(G)$ is a non-trivial matroid and 
$G \in {\cal G}_{\bowtie}$. Then $N$ is a component of $M_1(G)$ if and only if 
$N = M_1(A)$ for some component $A$ of $G$. 

\end{corollary}

\begin{remark}
\label{ReducToConnClean}

{\em 
Because of Corollary \ref{CmpM_1(G)}  we can reduce our problem of describing the classes of graphs having the same non-trivial bicircular matroid to the corresponding problem for the graphs in ${\cal CG}_{\bowtie}$, i.e. for connected cacti-graphs. 
}

\end{remark}

For this reasons from now on in the case of bicircular matroid we will assume that $G \in {\cal CG}_{\bowtie}$.

\vskip 1.5ex

From Corollary \ref{CmpM_1(G)} we have: 

\begin{Theorem} 
{\sc Graph description of connected matroid $M_1(G)$}
\label{M1Con}

\indent
The following are equivalent:
\\[1ex]
$(a1)$ $M_1(G)$ is a non-trivial matroid and $G \in {\cal CG}_{\bowtie}$ and 
\\[1ex]
$(a2)$ $M_1(G)$ is a connected matroid. 

\end{Theorem}

From Theorems \ref{maxk}, \ref{MkCon}, and \ref{M1Con} we have:

\begin{claim}
\label{2ConMkCon}

Let $G$ be a $2$-connected graph and $k \ge 1$. If $\Delta G \ge k$, 
then   $M_k(G)$ is a connected matroid.

\end{claim}

From Theorems \ref{CoreInGinfty}, \ref{maxk}, \ref{MkCon}, and
 \ref{M1Con} we also have:

\begin{claim}
\label{ConMGinfty}

Let $k \ge 1$. If $M_k(G)$ is a connected matroid, then 
$G \in {\cal G}_{\bowtie}$ and $\Delta G \ge k$.   

\end{claim} 

\begin{remark}
\label{FinRed}
{\em
The results in this section will allow us to reduce our problem $(WP)_k$ of describing the classes of graphs having the same $k$-circular matroid to the corresponding problem for cacti-graphs graphs, 
when $2 \le k \le \Delta G$, and for connected cacti-graphs,  
when $ k = 1 $. Equivalently, the theorems in this section allow us to reduce our problem $(WP)_k$ of finding all graph representations of a given $k$-circular matroid, $k \ge 1$,  to the corresponding problem for the class of connected $k$-circular matroids. 
}
\end{remark}

\subsection{Independent sets and bases  of  a $k$-circular 
\\
matroid}

\begin{claim}
\label{notreeind}

Let $G = (V, E,\phi )$ and $X \subseteq E$.
If $\Delta G\langle  X \rangle < k$ and $G\langle  X \rangle$ has no tree component, then $X \in {\cal I}_{k}(G)$.

\end{claim}

\bp (uses Claim \ref{deltaxleqdeltah})

Let $X\subseteq E$ and $\Delta G\langle  X \rangle < k$. If $ G\langle  X \rangle$ has no tree-component, then by Claim \ref{deltaxleqdeltah}, $\Delta G\langle  Z \rangle \leq \Delta G\langle  X \rangle < k$ for every $Z \subseteq X$. Hence $X$ contains no circuit of $ M_k(G)$ and therefore $X$ is an independent set of $M_k(G)$.  
\ep

\begin{claim}
\label{augmenting}

Let $G$ be graph and $k \geq 1$. Suppose that $I \in {\cal I}_{k}(G)$, $T$ is a tree component of 
$ G\langle  I \rangle $ or  $T$ is a vertex of $G$ not in 
$G\langle  I \rangle$, and $e$ is an edge in $E(G) \setminus I$ incident to at least one vertex of $T$. Then $I \cup e \in {\cal I}_{k}(G)$.   

\end{claim}

\bp (uses Claim \ref{cirstructure})

Suppose, to the contrary, that $I \cup e \notin {\cal I}_{k}(G)$. Then there 
exists a unique $C \in {\cal C}_{k}(G)$ such that 
$e \in C \subseteq I \cup e$. If $e$ has exactly one end-vertex in $T$, then the component of 
$ G\langle  C \rangle$ that contains $e$ will contain a leaf, contradicting Claim \ref{cirstructure}. If both ends of $e$ 
belong
to $T$ and $T'$ is the graph obtained from $T$ by adding edge $e$, then is connected and has at most one cycle.
Since the component $A$ of $ G\langle  C \rangle$ containing $e$ is a subgraph of $T'$, clearly $A$ also has at most one cycle contradicting Claim \ref{cirstructure} $(c2)$. 
 \ep

\begin{Theorem}  {\sc Graph structure of a base of  $M_{k}(G)$ in graph $G$} 
\label{basstructure}

Let $G$ be a graph and $k \geq 1$. 
Suppose that $M_k(G)$ is a connected matroid. 
Then the following are equivalent:
\\[1ex]
$(c1$) $B  \in {\cal B}_k(G)$ and
\\[1ex]
$(c2) $
$\Delta G \langle B \rangle = k - 1$, 
$V(G\langle  B \rangle) = V(G)$  $($i.e. $B$ spans 
$V(G)$$)$, 
and 
$\Delta A \geq 0$ 
for every component  $A$ of $ G\langle  B \rangle$
$($i.e. $ G\langle  B \rangle$ has no tree component$)$.

\end{Theorem}

\bp  (uses Claims \ref{deletingedge}, \ref{deltaxleqdeltaa},  \ref{deltaxleqdeltah}, \ref{geqkdep}, \ref{tree-components}, \ref{notreeind}, and \ref{augmenting})
\vskip 1ex
Since $M_k(G)$ is a connected matroid, by Theorem \ref{MkCon} for 
$ k \ge 2$ and Theorem \ref{M1Con} for $ k = 1$, matroid $M_k(G)$ is non-trivial and also graph $G$ is a cacti-graph. Therefore graph $G$ has no tree component. 
\\[1ex]
${\bf (p1)}$
First we prove that  $(c1) \Rightarrow (c2)$. 
Let $B \in {\cal B}_{k}(G)$.
\\[1ex]
\indent
${\bf (p1.1)}$
Our first step is to show that $V(G) = V(G\langle  B \rangle)$.
Suppose, to the contrary, that  there exists vertex $v$  in $V(G)\setminus V(G\langle  B \rangle)$. Since $G$ has no isolated vertices, there exists 
$e \in E(G) \setminus B$ such that $e$ is incident to $v$ in $G$.
Then by Lemma \ref{augmenting} (where $T$ is vertex $v$),
 $B \cup e\in {\cal I}_{k}(G)$. This   contradicts the maximality of $B$. Hence, if $B\in {\cal B}_{k}(G)$, then $V(G\langle  B \rangle) = V(G)$.  
\\[1ex]
\indent
${\bf (p1.2)}$
Our next step is to show that  $G\langle  B \rangle$ has no tree component. 

 \vskip 1ex
 
Suppose first  that $G\langle  B \rangle$ is a spanning tree. Since $M_k(G)$ is a non-trivial matroid, there exists $e \in E(G) \setminus B$.  By maximality of $B$, we have: $B \cup e$ is a dependent set of $M_{k}(G)$. Note that $\Delta_{G}(B \cup e) = 0$. 
Since $G\langle B \cup e \rangle $ is connected, by Claim \ref{deltaxleqdeltaa}. 
we have: $\Delta G\langle  X \rangle \leq \Delta G \langle B \cup e \rangle = 0 < k$ for every $X \subseteq B \cup e$. Therefore $B \cup e$ contains no $k$-circuit, a contradiction.  

 \vskip 1ex

Now suppose that $G\langle  B \rangle$ has more than one component and one of them is a tree component, say $T$. Since $G$ has no tree component, there exists $e \in E(G)\setminus B$ such that $e$ is incident to at least one vertex of $T$. Then by Lemma \ref{augmenting}, 
$B \cup e \in {\cal I}_{k}(G)$, a contradiction. 
Therefore 
$ G\langle  B \rangle$ has no tree component 
(i.e.
 $\Delta A \geq 0$ for every component $A$ of 
$G\langle  B \rangle$). 
\\[1ex]
\indent
${\bf (p1.3)}$
Now we show that $\Delta G\langle  B \rangle = k - 1$.
Suppose, to the contrary, that $\Delta G\langle  B \rangle \ne k-1$.
By Claim \ref{geqkdep},  $\Delta G\langle  B \rangle < k$. 
 Therefore $\Delta G\langle  B \rangle < k - 1$. 
 Since $M_k(G)$ is a non-trivial matroid, there exists edge $e$ in $E(G)\setminus B$. By Claim \ref{deletingedge}, $\Delta G\langle B \cup e \rangle \leq k-1 < k$. 
 Since $V(G\langle  B \rangle) = V(G)$ and $G$ has no tree component,
 clearly $G\langle  B \rangle$ has no tree component as well.
 Now by Claim \ref{notreeind},  $B \cup e \in {\cal I}_{k}(G)$. Therefore $B$ is not a maximal independent set of $M_k(G)$, and so $B$ is not a base of 
 $M_k(G)$, a contradiction. Thus, $\Delta G\langle  B \rangle \geq k - 1$,  that contradicts our assumption that $\Delta G\langle  B \rangle < k - 1$.
\\[1.5ex]
${\bf (p2)}$ Finally, we prove
$(c1) \Leftarrow (c2)$.  
Since no component of $G\langle  B \rangle$ is a tree,
by Claim \ref{deltaxleqdeltah},
 $\Delta G\langle  X \rangle \leq \Delta G\langle  B \rangle = k - 1 < k$ for every $X\subseteq B$. Therefore $B$ contains no circuit of 
 $M_k(G)$, and so $B \in {\cal I}_{k}(G)$. Since $V( G\langle  B \rangle = V(B)$, we have: $\Delta G\langle B\cup e \rangle > \Delta G\langle  B \rangle = k - 1$ for every $e\in E(G)\setminus B$. Then, 
 by Claim \ref{geqkdep}, $B \cup e \in {\cal D}_{k}(G)$ for every $e \in E(G) \setminus B$. Therefore $B$ is a maximal independent set of $M_k(G)$. 
 \ep

 \vskip 2ex

Let 
$\rho _k(G)$ and $\rho  ^*_k(G)$ denote the rank and the corank of matroid $M_k(G)$, respectively, and so 
$\rho _k(G) + \rho ^*_k(G) = |E(G)|$.
Obviously, if ${\cal C}_{k}(G) = \emptyset $, then 
$\rho _k(G) = |E(G)|$ and $\rho ^*_k(G)  = 0$. We also know that
if $G$ has no isolated vertices, then  $\rho _0(G) = |V(G)| - cmp(G)$, 
where $cmp(G)$ is the number of components of $G$.

\begin{corollary}
\label{grank}

Let $G$ be a graph and $k \geq 1$. Suppose that $M_k(G)$ is a connected matroid. 
\vskip 0.6ex
\noindent 
Then 
$\rho _k(G) = |V(G)|  - 1 + k $ and  $\rho ^*_k(G) = |E(G)| - |V(G)| + 1 - k$.
 
\end{corollary}

 \bp (uses Claim \ref{basstructure})
 
By Claim \ref{basstructure} $(c2)$, if $B$ is a base of $M_k(G)$, then

\vskip 0.4ex
\noindent 
$|B| - |V (G  \langle B \rangle)| =  \Delta G \langle B \rangle = k - 1$ and $V(G\langle  B\rangle) = V(G)$. Therefore
\vskip 0.4ex
\noindent 
$\rho _k(G) = |B| = |V (G  \langle B \rangle)|  - 1 + k =
|V(G)| - 1 + k$ and
$\rho ^*_k(G) = |E(G)| - |V(G)| + 1 - k$.
  \ep

\subsection{The core of a base of a $k$-circular matroid }
\label{core}

\indent

In this Subsection we define the notion and describe some properties of the core of a base of $M_k(G)$ as a special subgraph of $G$. This notion will be essential for establishing some properties of the cocircuits of $M_k(G)$ in Subsection  \ref{cocircuits}. 

\vskip 1.5ex

The definition below uses the notion of the core $ [G ]$ of a graph $G$ (see Definition \ref{NonConGraphCore}).

\begin{definition} {\sc The core of a base of $M_k$(G)}.
\label{CoreOfB}

{\em 
Suppose that $k\ge 2$, $M_k(G)$ is a connected matroid, and $B$ is a base of $M_k(G)$. 
\\
Then the core $ [G \langle B \rangle]$ of graph $G$ is also called {\em the core of $B$ in $G$}. 
}
\end{definition}

By Theorem \ref{basstructure}, if $k \ge 2$ and $B$ is a base of  a connected matroid $M_k(G)$, then the core 
$ [G \langle B \rangle] $ of graph $ G \langle B \rangle $ is defined and if $k = 1$, then the core of $B$ in $G$ is not defined. 

\vskip 1.5ex

From Theorems \ref{CoreInGinfty} $(a3)$ and \ref{cirstructure} we have:

\begin{claim}
\label{DescCoreOfB} 

Let $k \ge 2$ and $B \in {\cal B}_k(G)$. Suppose that $M_k(G)$ is a  connected matroid. 
\vskip 0.6ex
Then 
\\[1ex]
$(c1)$ $[G \langle B \rangle]$ is the unique subgraph of 
$ G \langle B \rangle $ belonging to ${\cal G}_{\bowtie} $ such that 
\\[0.6ex]
$\Delta  [G \langle B \rangle]  = \Delta G \langle B \rangle = k -1 $ and 
\\[1ex]
$(c2)$ the edge set of $ [G \langle B \rangle] $ is the unique subset of $B$ that is a  $($$k$-$1$$)$-circuit of $G$.   
 \end{claim}

\begin{claim}
\label{k>1BasesPrevCirc} 

Let $k \ge 2$ and $C \in {\cal C}_{k-1}(G)$. Suppose that $M_k(G)$ is a  connected matroid. 
\vskip 0.6ex
Then there exists $B \in {\cal B}_k(G)$ such that $G \langle C \rangle = [G \langle B \rangle] $.

\end{claim}

\bp (uses Theorem \ref{basstructure})
\vskip 0.6ex
Let $Q \in Cmp(G)$. If $C \cap E(Q) \ne \emptyset$, then let 
$Q'$ be a minimal connected spanning subgraph of $Q$ containing $C$. If 
$C \cap E(Q) = \emptyset$, then let $Q'$ be a minimal connected unicyclic subgraph of $Q$. Clearly, $H = \cup \{Q': Q \in Cmp(G) \} $ satisfies condition 
$(c2)$ of Theorem \ref{basstructure}. Therefore $B = E(H) \in {\cal B}_k(G)$ and $G \langle C \rangle = [G \langle B \rangle] $. 
\ep

\vskip 1.5ex
From Claims \ref{DescCoreOfB} and \ref{k>1BasesPrevCirc} we have: 

\begin{claim}
\label{k>1Core=PrevCirc} 

Let $k \ge 2$. Suppose that $M_k(G)$ is a  connected matroid. Then 

$ {\cal C}_{k-1}(G)  =  \{ [G \langle B \rangle]: B \in {\cal B}_k(G)\}$.

\end{claim}

The claim below uses the notion of the kernel $\lfloor G \rfloor$ of a graph 
$G$ (see Definition \ref{ConGraphCore}).

\begin{claim}
\label{k=1BasesPrevCirc}

Let $k = 1$. Suppose that $M_1(G)$ is a  connected matroid. Then 
\\[1ex]
$(c1)$ if $B \in {\cal B}_1(G)$, then $\lfloor A \rfloor$ is a cycle for every 
$A \in Cmp (G \langle B \rangle)$ and 
\\[1ex]
$(c2)$ if $C$ is the edge set of a cycle of $G$, then there exists 
$B \in {\cal B}_1(G)$ and $A \in Cmp (G \langle B \rangle)$ such that 
$G \langle C \rangle = \lfloor A \rfloor $.

\end{claim}

\bp (uses Theorem \ref{basstructure})
\vskip 0.6ex
First we prove $(c1)$.By Theorem \ref{basstructure}, 
$\Delta G \langle B \rangle = 0$ and $\Delta A \geq 0$ for every $A \in Cmp (G \langle B \rangle)$. Now 
$G \langle B \rangle = \cup \{A: A \in Cmp (G \langle B \rangle) \}$, and so
$\Delta G \langle B \rangle = \sum  \{\Delta A : A \in Cmp (G \langle B \rangle)\}$. 
Therefore $\Delta A = 0$ for every $A \in Cmp (G \langle B \rangle)$. Then every component $A$ of $G \langle B \rangle$. Thus $\lfloor A \rfloor$ is a cycle for every $A \in Cmp (G \langle B \rangle)$.

 \vskip 1ex
 
Now we prove $(c2)$. Let $C$ be a cycle of $G$. Let $A \in Cmp(G)$. If $A$ contains $C$, then let $A'$ be a minimal connected spanning subgraph of $A$ containing $C$. If $A$ does not contain $C$, then let $A'$ be a minimal connected unicyclic spanning subgraph of $A$ containing $C$. Clearly, 
$H = \cup \{A': A \in Cmp(G) \} $ satisfies condition $(c2)$ of Theorem \ref{basstructure}. Therefore $B = E(H) \in {\cal B}_k(G)$ and 
$G \langle C \rangle = \lfloor A \rfloor $ for some $A \in Cmp(G)$. 
\ep

\subsection{Cocircuits of a $k$-circular matroid}
\label{cocircuits}

\indent

In this Section  we intend to describe the  cocircuits of the 
$k$-circular matroid $M_k$ of a graph $G$ as some special edge subsets of $G$. 
For example, we know that $K$ is a cocircuit of $M_0(G)$ if and only if $K$ is a minimal edge cut in $G$.

 \vskip 1ex
 
By Claim \ref{FundCircuit} $(c3)$,
for every cocircuit $K$ of a matroid $M$ there exists 
a base $B$ of $M$ and $e \in B$ such that $K $ is the fundamental  
cocircuit of base $B$ in $M$ rooted at $e$ (or the same, a $(B,e)$-cocircuit in $M$), i.e. $K = C^*(e,B)$ for some $B \in {\cal B}(M)$ and $e \in B$. Put $K(e,B) = C^*(e,B)$.

Thus, we can (and will) describe all cocircuits of matroid $M_k(G)$ by describing the $(B,e)$-cocircuits $K(e,B)$ of 
$M_k(G)$ for all pairs $(B,e)$, where $B$ is a base of $M_k(G)$ and $e \in B$. 

 \vskip 1ex
 
We will distinguish between three possible types of $(B,e)$-cocircuits $K(e,B)$ depending on the structure of component $A$ of $ G \langle B \rangle $  containing edge $e$ and on the position of edge $e$ in $A$.
\begin{definition} 
\label{cocircuits-type}
Let $B \in {\cal B}_k(G)$, $e \in B$, and $k \ge 1$. Then
\\[0.7ex]
$(t1)$ $K(e,B)$ is a {\em $(B,e)$-cocircuit  in $M_k(G)$ of type 1},
if $e \notin E \lfloor A \rfloor $, where $A$ is a component of 
$ G \langle B \rangle$ containing edge $e$,
\\[0.7ex]
$(t2)$ $K(e,B)$ is a {\em $(B,e)$-cocircuit  in $M_k(G)$ of type 2},
if $e \in E \lfloor A \rfloor $, where $A$ is a unicyclic component of 
$ G \langle B \rangle$ containing edge $e$, and
 \\[0.7ex]
$(t3)$ $K(e,B)$ is a {\em $(B,e)$-cocircuit  in $M_k(G)$ of type 3},
if $e \in E [ A ] $, where $A$ is a component of 
$ G \langle B \rangle$ that has at least two cycles and contains edge $e$.
 \end{definition}

 \begin{Theorem} {\sc Graph description of rooted cocircuits of type 1}
 \label{type1}

Let $G$ be a graph and $k \ge 1$. Suppose that  
\\[1ex] 
$(a1)$ $M_k(G)$ is a connected matroid, $B \in {\cal B}_k(G)$, and $e \in B$ and  
\\[1ex] 
$(a2)$  $K(e,B)$ is a $(B,e)$-cocircuit  in $M_k(G)$ of type 1 $($i.e. edge $e \notin E \lfloor A \rfloor $, where $A$ is a component of 
$ G \langle B \rangle$ containing edge $e$$)$.  

\vskip 1ex
Then exactly one of the two components of $A \setminus e$ is a tree $T$ 
and $K(e,B) = K'(e,B) \cup e$, where $K'(e,B)$ is the set of edges in  
$E \setminus B$ having at least one end-vertex in $V(T)$.

 \end{Theorem} 
 
 \bp  (uses Theorem  \ref{basstructure} and Claims \ref{FundCircuit} and 
 \ref{A-e}) 
  \\[1ex]
${\bf (p1)}$
Since $e \notin E \lfloor A \rfloor $, by Claim \ref{A-e} $(c1)$, exactly one of the two components of $A \setminus e$ is a tree. Let $T$ be the tree component and $D$ the non-tree component of 
$A \setminus e$.   
 \\[1.5ex]
${\bf (p2)}$
We prove that $K(e,B) = K'(e,B) \cup e$.
 By Claim \ref{FundCircuit} $(c2)$, it is sufficient to show that 
 $u \in  K'(e,B)$ if and only if $u \in E \setminus B$ and 
 $B_u = (B \setminus e) \cup u$ satisfies condition $(c2)$ of Theorem  \ref{basstructure}.
\\[1ex]
${\bf (p2.1)}$ First, suppose that  $u \in  K'(e,B)$.
 
If $u$ has both ends in $V(T)$, then both $D$ and $T \cup u$ have cycles, and therefore  $B_u$ satisfies condition $(c2)$ of Theorem \ref{basstructure}.
  
  If $u$ has one end in $T$ and the other end in $D$, then
 $ G \langle (E(A)\setminus e) \cup u \rangle $ is not a tree and $B_u$  satisfies condition $(c2)$ of Theorem \ref{basstructure}. 
  
 If $u$ has exactly one end in $T$ and  the other end in a component 
 $F $ of $G \langle B \rangle $ distinct from $A$,
 then $ G \langle E(F) \cup E(T) \cup u \rangle$ is not a tree, 
 and therefore $B_u$ satisfies condition $(c2)$ of Theorem \ref{basstructure}.  
 \\[1ex]
${\bf (p2.2)}$
  Now suppose that $u \in  E \setminus K(e,B)$.
Then $T$ is a tree-component of 
$G\langle B\setminus e \cup  u \rangle $, and so $B_u$ 
does not satisfies condition $(c2)$ of Claim \ref{basstructure}.  
 \ep

\begin{Theorem} {\sc Graph description of rooted cocircuits of type 2}
 \label{type2}

Let $G$ be a graph and $k \ge 1$. Suppose that  
\\[1ex] 
$(a1)$ $M_k(G)$ is a connected matroid, $B \in {\cal B}_k(G)$, and $e \in B$ and  
\\[1ex] 
$(a2)$ $K(e,B)$ is a $(B,e)$-cocircuit  in $M_k(G)$ of type 2, $($i.e. 
$e \in E \lfloor A \rfloor $, where $A$ is a unicyclic component of 
$ G \langle B \rangle$ containing edge $e$$)$.    

\vskip 1ex
Then $A \setminus e$ is a tree and $K(e,B) = K'(e,B) \cup e$, where $K'(e,B)$ is the set of edges in $E \setminus B$ having at least one end-vertex in 
$V(A \setminus e)$.

 \end{Theorem} 
 
 \bp  (uses  Theorem  \ref{basstructure} and Claim \ref{FundCircuit})

By assumption $(a2)$ of our theorem, $e$ is an edge of a unique cycle in $A$.
If $u$ has both ends in $V(A)$, then 
$\langle (E(A) \setminus e) \cup u \rangle $ is not a tree and therefore 
$(B \setminus e) \cup u$ satisfies condition $(c2)$ of Theorem \ref{basstructure}. 
 
 If $u$ has exactly one end in $V(A)$ and the other end in a component $D$ of 
 $G \langle B \setminus e \rangle$, then   $ G \langle E(D) \cup E(A\setminus e) \cup u \rangle $ is not a tree and $(B\setminus e) \cup u$ satisfies condition 
 $(c2)$ of Theorem  \ref{basstructure}.   
\\
\indent
Now, if $u $ has no ends in $E(A)$ or, equivalently, $u \notin K'(e,B)$, 
then $A\setminus e$ is a tree component of 
$G \langle (B\setminus e) \cup u \rangle$, and so $(B\setminus e) \cup u$ does not satisfies condition $(c2)$ of Theorem \ref{basstructure}. Therefore 
$(B\setminus e) \cup u \notin {\cal B}_{k}(G)$. 

Hence by  Claim \ref{FundCircuit} $(c2)$, 
\\
$u \in K'(e,B)\cup e = K(e,B) 
\Leftrightarrow (B\setminus e) \cup u \in {\cal B}_{k}(G) 
\Leftrightarrow u \in K(e,B)$.   
\ep

\begin{Theorem} {\sc Graph description of rooted cocircuits of type 3}
 \label{type3}

Let $G$ be a graph and $k \ge 1$. Suppose that  
\\[1ex] 
$(a1)$ $M_k(G)$ is a connected matroid, $B \in {\cal B}_k(G)$, and $e \in B$ and  
\\[1ex] 
$(a2)$ $K(e,B)$ is a $(B,e)$-cocircuit  in $M_k(G)$ of type 3, $($i.e. 
$e \in E [ A ] $, where $A$ is a component of $ G \langle B \rangle$ that has at least two cycles and contains edge $e$$)$.    

\vskip 1ex
Then $K(e,B) = (E \setminus B) \cup e$.  
 
 \end{Theorem}

 \bp  (uses  Theorem \ref{basstructure} and  Claim \ref{A-e})

Since $e \in E [ A ] $, by Claim \ref{A-e} $(c3)$, every component of $A \setminus e$ contains a cycle. Therefore for every $u \in E\setminus B$ we have: 
$(B\setminus e) \cup u$ satisfies condition $(c2)$ of Theorem \ref{basstructure}. Hence $u \in K(e,B)$ if and only if $u \in (E\setminus B) \cup e$. 
\ep

\end{document}